\documentclass[dvips,12pt]{article}
\usepackage{epsfig,graphicx}

\usepackage{latexsym}
\usepackage[psamsfonts]{amssymb}
\usepackage{amsfonts}
\usepackage{amssymb}
\usepackage{amsmath}

\def\r{\mathbb R}
\def\e{\mathbb E}
\def\p{\mathbb P}

\def\c{\mathbb C}

\def\CP{{\cal P}}

\def\CR{{\cal R}}

\def\re{{\mathfrak {Re}\,}}
\def\im{{\mathfrak{Im}\,}}

\def\defn{\stackrel{\rm def}{=}}

\newtheorem{theorem}{{\bf Theorem}}%[section]
\newtheorem{lemma}{{\bf Lemma}}%[section]
\date{}

\begin{document}
\newcommand{\proof}{\noindent {\sl Proof.  }}%[section]
\newcommand{\sproof}{\noindent {\sl Scetch of the proof.  }}%[section]
\newcommand{\rem}{\noindent {\sl Remark.  }}%[section]
\title{On  Stable Pareto Laws \\ in a
Hierarchical Model of Economy\footnote{This work
is supported by
DFG 436 RUS 113/779/0-1.}}
\author{A.~M.~Chebotarev\\ {\small \it Quantum Statistics Department,
Moscow State
University},\\
{\small \it Moscow 119899, Russia}}
%e-mail: $\rm a_{-}m_{-}$chebotarev\@ newmail.ru
\maketitle

\begin{abstract}
This study considers a model of the
income distribution of agents whose
pairwise interaction is asymmetric and
price-invariant. Asymmetric
transactions are typical for
chain-trading groups who arrange their
business such that commodities move
from senior to junior partners and
money moves in the opposite direction.
The price-invariance of transactions
means that the probability of a
pairwise interaction is a function of
the ratio of incomes, which is
independent of the  price scale or
absolute income level. These two
features characterize  the hierarchical
model. The income distribution in this
class of models is a well-defined
double-Pareto function, which possesses
Pareto tails for the upper and lower
incomes. For  gross and net upper
incomes, the model predicts definite
values of the Pareto exponents, $a_{\rm
gross}$ and $a_{\rm net}$, which are
stable with respect to quantitative
variation of the pair-interaction. The
Pareto exponents are also stable with
respect to the choice of a demand
function within two classes of
status-dependent behavior of agents:
linear demand ($a_{\rm gross}=1$,
$a_{\rm net}=2$) and unlimited slowly
varying demand  ($a_{\rm gross}=a_{\rm
net}=1$). For the sigmoidal demand that
describes limited returns, $a_{\rm
gross}=a_{\rm net}=1+\alpha$, with some
$\alpha>0$ satisfying a transcendental
equation. The low-income distribution
may be singular or vanishing in the
neighborhood of the minimal income; in
any case,  it is $L_1$-integrable and
its Pareto exponent is given
explicitly.

The theory used in the present study is based on a simple
balance equation and new results from
multiplicative Markov chains and
exponential moments of random
geometric progressions.
\end{abstract}

\section{Introduction}
\subsection{A brief historical survey}
The orthodox interpretation of
statistical distributions in economics
assumes a certain
stochastic process whose invariant
probability measure characterizes the
given distribution. Champernowne
\cite{Ch53} made a remarkable
attempt to model the Paretian
distribution, and obtained several
discrete versions of the distribution
as  stationary probability measures for
a class of Markov chains on half-line.
Starting from an empirical analysis of
UK tax records for the period 1951--52,
Champernowne reconstructed a $(12\times
15)$-approximation of the stochastic
matrix of the Markov chain, and noted
a definite ``tendency for the lowest
income to shift upwards'' \cite{Ch53}.
In a sense, his basic
assumptions (asymmetry and homogeneity
of the transition probability) are inherited
by the model described in the present study.
To ensure the probabilistic
interpretation and to obtain a closed
equation for the stationary
distribution, Champernowne introduced a
boundary condition for the occupation
number of the lowest income range. This
fare for mathematical correctness was
not reasonably explained in economic
terms, but a similar model of
multiplicative stochastic process with
``reset'' events (i.e., with stochastic
renewal of the process at the boundary) \cite{MaZa99}
was successfully applied in \cite{NS04}
for a Monte Carlo simulation of
variations in income distribution
within Japan. In the present paper,
only the continuity
assumption is imposed at the boundary (at infinity
and at  the point corresponding to the
minimal income) that defines  a stationary density distribution.

To fix notation, recall
that the cumulative probability
distribution ${\rm Pr}(x\ge s)$
of a random variable $x$ has the
Pareto tail with exponent $a$ if
${\rm Pr}(x\ge s)=O(s^{-a})$, $a>0$, for large
$s$\footnote{In fact, the conventional
terminology is not completely correct:
Zipf constructed his plots with the probability
axis $y$, while Pareto used $x$ for the
same purpose \cite{New05}.}.

A number of empirical facts in
economics  concerning the Pareto
distribution  can be qualitatively
explained using generalizations of the
Yule stochastic process \cite{Yul37},
which was proposed in 1937 as a simple
dynamical model of  biological taxa. In
view of a certain similarity between
the present model and Yule's model, it
is worth explaining the basics of
Yule's model  in a macroeconomic
context. Suppose that the starting
capital equals $1$ and a nonrandom
consumer demand provides the business
community  with fixed $1+m$ units of
income per unit of time. Each step of
the time-evolution consists of two
events: (1) a new agent obtains  the
starting capital and joins the business
community; (2) the remaining $m$ units
of income are randomly distributed
among the existing agents according to
their status-demand, for example, with
partial probabilities proportional to
the existing distribution of  capital
(linear demand). The construction of  a
stationary solution in a biological
context is discussed in \cite{New05}.
The cumulative stationary distribution
has the power tail with the following
exponent:
$$
a=1+\alpha, \quad \alpha= \frac{1}{m}.
$$
If the business society only rarely recruits new
agents, $\alpha$ is
small and $a\approx 1$\footnote{Let $n$
be the number of time units per year.
Then, $n\,m\approx n\,(1+m)=\sigma$ is
the size of the market and $\nu=n$ is
the intensity of the creation of vacancies.
If $a\ge 1$, it is considered that  $\alpha=a-1$ is the
{\it horizon} of the power tail,
as the heaviest tail with $a=1$
normally  overshadows the remaining
components of a sampling. In these
terms, Yule's law states: ``{\it Horizon
$\times$ market size $\thicksim$
intensity of vacancies' creation}''.}.

In addition to linear demand
($\ell$-demand), the class
of slowly varying functions is considered. From
an empirical viewpoint, in the linear
case the size of a business is far
from market limitations or
the constraints of resources. The slowly varying
demand ($u$-demand) carries
features of utility and sigmoidal
functions that describe the diminishing
returns of investments and  the effects
of market  or production constraints.

Empirical studies of the top income
distribution show that the Pareto
exponent $a$ of the cumulative
distribution normally belongs in the
interval $[1,\,3]$. There is evidence
for both stability
\cite{AoSouF03}--\cite{A04} and
instability \cite{AtSal02, At03} for
this indicator over time. In the latter
two reviews, the authors discuss the
relationship between tax legislation
and the distribution of personal
income. The reviews present variation
of the Pareto exponent in the range
[1.3, 3.5], based on tax records from
Great Britain and the Netherlands for
the years 1907--1999, and records from
the US, Canada, and France over shorter
periods. The detailed records enable us
to monitor details of income
distribution for the top 0.05\% to 10\%
of earners. The Pareto index $a=1$ for
the net income of firms and
corporations  has been observed by many
authors (see \cite{FuGuA05, MKTT,
KII}), including the author of the
present paper \cite{Ch03, Ch05}. These
tails are the least informative in the
sense that they overshadow the possibly
important details of less heavy tails
that contribute to the general picture.

%\psdraft

\vskip5mm
\begin{figure}[h]
\begin{minipage}{0.5\linewidth}
%\centering\epsfig{figure= S0x0p0,width=\linewidth}
\centering\epsfig{figure= 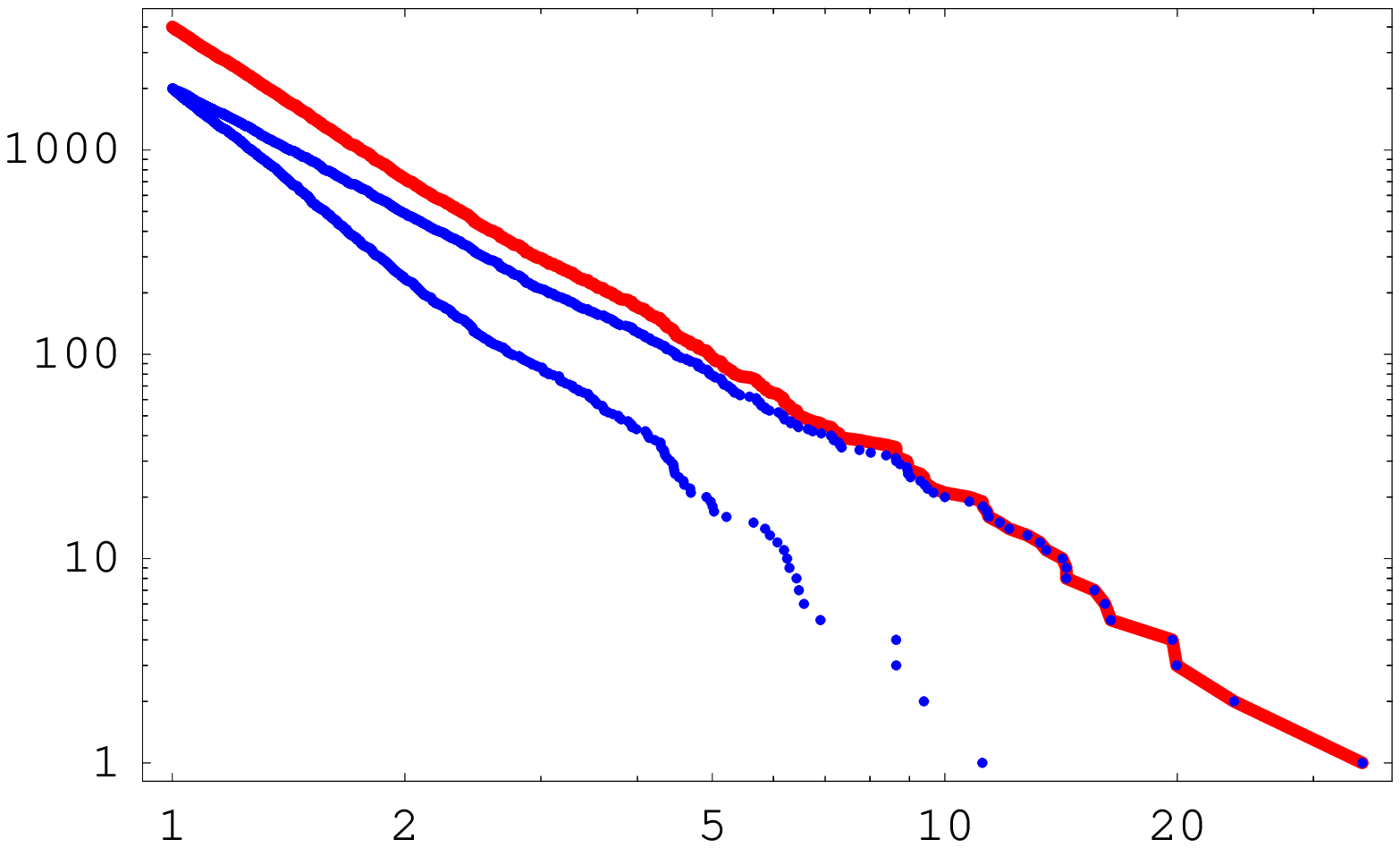,width=\linewidth}
\end{minipage}\hfill
\begin{minipage}{0.5\linewidth}
%\centering\epsfig{figure= S0x0p0around0.EPS,width=\linewidth}
\centering\epsfig{figure= 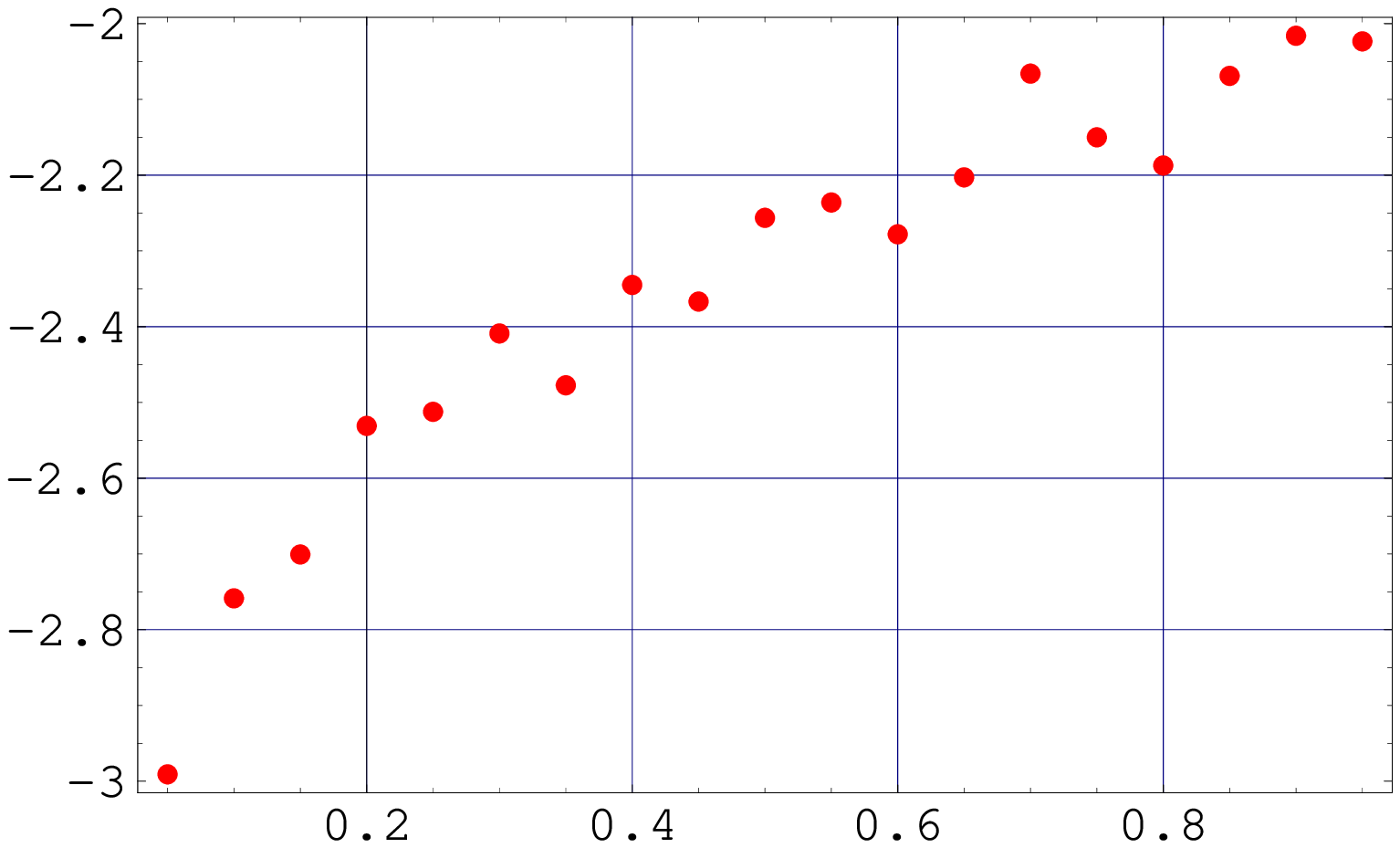,width=\linewidth}
\end{minipage}\hfill
\caption{\small The two lower curves in the left-hand
figure display the cumulative
probabilities
defined by 2000 Pareto-distributed
random points with exponents $a=1$ and
$a=2$, respectively, while the upper curve
shows their 1:1 mixture: $a_{\rm
mix}\approx 2.3$. This illustrates the shadowing effect
created by the heaviest tail.  The point-plot on
the right shows $a_{\rm mix}(x)$ as a
function of the share $x$ of the
heaviest
tail varying from 5\% to 95\% over the
total sampling of  4000 random
points.}
\end{figure}

A solvable model for the Pareto
income distribution was recently proposed  by
Reed \cite{R03}. In fact, Reed's
approach is based on the  Poincare
ergodic conjecture, which states that
the time-average coincides with the
average over the  ensemble.  The
ensemble considered by both Reed and Gabaix
\cite{Gab99} consists of
noninteracting agents whose income is
described by geometric Brownian
motion and whose life-time is an
exponentially distributed random
variable. By averaging the individual
income distribution over the random
life-time, the
double-Pareto distribution
\cite{R03} is derived\footnote{A correspondent
informed me that in 1949, E.~Fermi described a similar mechanism
in an astrophysical context, but I have been unable
to confirm this.}. According to the
Poincare conjecture, the obtained
distribution describes the income
distribution over the  ensemble of
ideal (i.e., noninteracting) agents.

The present paper considers an ensemble
of asymmetrically interacting agents.
The balance equation that describes the
stationary distribution of income is
derived using arguments that are
similar to those used by Bouchaud
\cite{Bush01} for symmetric
pair-interactions. The present equation
enables the calculation of the pair of
Pareto exponents $a$ and $d>-1$ for the
top and bottom income tails. The value
of the Pareto exponent $a_{\rm net}=1$,
derived below for hierarchical models
with $u$-demand, is normal for the top
revenues of firms and corporations in
Japan (see \cite{NS04, AoSouF03}) and
is in agreement with estimates  for
Europe, Britain, and the USA (see
\cite{A04, At03, R03, KK03}). Indirect
reconstruction of the gross income
distribution for the top 5\% of
householders, based on sale statistics
for the national car market and  direct
estimates of the capital distribution
among the 200 largest banks in Russia,
yields the same value
\cite{Ch03, Ch05}. A possible
qualitative explanation for this rather
general phenomena  is the existence of
a dominating corporate subculture of
profit-making communities. The  aim of
this paper is to consider the extreme
case of a corporate community: the
hierarchical structure. A  sketch of
this approach was published in
\cite{Ch04}.

\subsection{Assumptions and results}
In the  business world, the
hierarchical environment maintains and
permanently reproduces the specific
circulation of money and commodities;
corporations, network-marketing and
chain-trading structures, shadow
economy groups, etc., aim to organize
supply and demand such that
commodities and money  move in
opposite directions.    Within the
hierarchy,  an agent consumes  the commodities
distributed by senior partners, while at
the same time, he distributes them to
junior agents \cite{Bog99}.

Clearly, the above picture is not
strictly true for agents with low
income  who are receiving money  (e.g., low
wages or pensions, welfare subsidies or
grants) from the state.
To simplify mathematical
considerations,
taxation and
the redistribution policy of  the state is ignored,
and a  toy model is considered
that describes the stable distribution of
unilateral money flow in  a
hierarchical structure under
the assumption that the intensity
$\CR(s,x)$ of transactions between
the pair of agents is an
asymmetric and price-invariant
function of incomes $s$ and
$x$ of these agents; that is,
$\CR(s,x)=R(s/x)\ge0$ if
$s\ge x$, and $R(s/x)=0$ if
$s<x$.  These assumptions
characterize this class of hierarchical
models. We
also assume that income is
(mathematically) unlimited and that the minimal income $s_0$
is positive, i.e., $s\in[s_0,\infty)$,
$s_0>0$.

The money flow controlled by an agent
with net income $s$ consists of the
welfare payment  $\CP(s_0,s)$ (cash
benefits, social grants, etc.) and
the contribution $B(n,s)$ from
subordinate agents. Here, $n=n(s)\in
L_1(s_0,\infty)$ is the density of
agents with net income $s$, and $B(n,s)$
is a function of $s$ and a functional
of $n$. The hierarchy assumption means
that $B(n,s)$ depends on the values of
$n(x)$ restricted to $x\in [s_0,s]$,
and $B(n,s)\to 0$ as $s\to s_0$. The welfare $\CP(s_0,s)$
decreases and vanishes for all
sufficiently large $s$. The conjectured
price-invariance of the model implies
that $\CP(s_0,s)=P(s/s_0)$.

To maintain the social status,
an agent with net income $s$ pays out
$H=H(n,s)$ to   senior partners. We
assume that  the hierarchy is free of
competition; in other words, we
conjecture that   (at least in the
neighborhood of the stable
distribution)
variations in the
number of senior agents does not affect
consumer demand and price levels.
This assumption implies that the
individual
consumption costs
$H=H(s)$  are independent of
$n$. Thus, the {\it net income} is equal
to
\begin{align}
\label{balance}
s=P(s/s_0)+B(n,s)-H(s),
\end{align}
and the sum
$g(s)=P(s/s_0)+B(n,s)=s+H(s)$ is
the   {\it gross income}. Suppose that the consumption costs
$H(s)$ decrease as $s\downarrow s_0$. Hence, the
gross income $g(s)=s+H(s)$ is an increasing
function of $s$ and the  equation
\begin{align}
\label{eqS0}
s_0:\; s_0+H(s_0)=P(1)
\end{align}
has a unique
solution. Furthermore,
(\ref{eqS0})
is regarded as an analytical definition of the minimal income $s_0$,
while (\ref{balance}) is considered as an
equation with respect to  $n(s)$.

In Section 2, $B$ and $H$ are represented
in terms of
the pair-interaction intensity $R(x)$ and
the demand function $\sigma(s)$.
The main statement
of this section is a conditional
proposition: if  the solution $n(s)$ of
equation (\ref{balance}) can be
represented in the form
\begin{align}
\label{asympt}
n(s)&=x^b\,(1-x)^d\,
\rho(x)/\sigma(x^{-1}),\quad x={s_0}/{s}\in(0,1]
\end{align}
where $\rho$  is a
continuous, strictly positive, and uniformly bounded
function on $[0,1]$, then
\begin{align}
\label{d}
d= -1+\frac{R(1)}{C(1)}> -1,
\quad C(1)=1-P'(1)+2R_0
>0,\quad R_0\defn \int_1^\infty
R(x)\,dx.
\end{align}
In addition, it is demonstrated that $a_{\rm
gross}=1$, $a_{\rm net}=2$ for the
linear demand, $a_{\rm
gross}=a_{\rm net}=1$ for any slowly
varying demand, and $a_{\rm
gross}=a_{\rm net}=1+\alpha$ for the sigmoidal demand, where
$\alpha$ is a unique solution of the
transcendental equation
$$
\alpha:\;\frac{1}{R_0}\int_1^\infty
s^\alpha R(s)\,ds=1+\frac{1}{R_0S_0},\quad S_0=\lim_{s\to\infty}
\sigma(s).
$$

Section 3 considers  the results of
numerical simulation of
$\rho$ that illustrate various
opportunities of the model.
The existence of the income
distribution $n(s)$ in the  form (\ref{asympt}) is
proved  in Section 4. $\rho$ is first represented in
terms of   the  pair of multiplicative
Markov chains, $\{\xi_n\}$ and
$\{\eta_n\}$,  that describe the function
$\rho(s)$ either as $s\to \infty$ or
$s\to s_0$.
In Section 5, possible
economical issues  are discussed.
Section 6 consists of a mathematical appendix
that contains four basic lemmas on
multiplicative Markov chains. It is anticipated
that mathematical facts considered in
the last
section will be useful for the
readers who use multiplicative Markov processes in econophysics
and mathematical finance.

\section{Balance equation and continuity conditions}
\subsection{Balance equation}
Suppose that the
money flow from agents
with income $x$ to $n(s)$ agents with
higher income $s$ is proportional to
the density
$n(x)$, consumer demand
$\sigma(x)$, and  the
pair-potential $R(s/x)$ of interactions.
As this income is shared among $n(s)$ senior agents, the
personal inflow functional $B(n,s)$ can be written as follows:
$$
B(n,s)=\frac{1}{n(s)}\int_{s_0}^s
\sigma(x)\,R(s/x)\,n(x)\,dx.
$$
The given rational agent behaves in
exactly the same way toward the senior
partners as his junior partners:  the
cost of his demand to agents with
higher income $x$ is proportional to
the intensity of interactions  and to
his demand function $\sigma(s)$:
$$
H(s)=\sigma(s)\, \int_s^\infty
R(x/s)\,dx=s\,\sigma(s)\,R_0,\quad R_0\defn
\int_1^\infty R(y)\,dy.
$$
Taking into account the welfare income $P$,
the balance equation (\ref{balance}) is rewritten as follows:
\begin{align}
\label{income}
s=P(s/s_0)+\frac{1}{n(s)}\int_{s_0}^s
\sigma(x)\,R(s/x)\,n(x)\,dx-
\sigma(s)\, \int_s^\infty R(x/s)\,dx.
\end{align}
This relation is considered
as an equation for the stationary
distribution $n(s)$ of the net income $s$.
The equation
$
g=s\,(1+\sigma(s)\,R_0)
$
is an obvious statement on the relation
between the {\it gross} and {\it net income}, $g$ and
$s$, which readily follows from
(\ref{income}).

For further use of the model, it is convenient
to rewrite (\ref{income}) as  a linear
homogeneous integral equation:
\begin{align}
\label{stat}
\int_{s_0}^s
\sigma(x)\,R(s/x)\,n(x)\,dx
=\bigl(s-P(s)+s\,\sigma(s) R_0\bigr)\,n(s).
\end{align}
The choice  of a smooth monotone demand
$\sigma(s)\le s$ is restricted to the
following three classes:
\begin{itemize}
\item[($\ell$):] $\sigma(s)=s$ for
all $s>1$, \;\; ($s$): $\sigma(s)=s$ for $s=O(1)$;
$\sigma(s)\to S_0 \text{\;
as\;  } s\to \infty$,
\item[($u$):] $\sigma(s)=s$ for $s=O(1)$;
$\frac{\sigma(s t)}{\sigma( s)}\to 1,
\;\sigma(s)\to \infty \text{\;
as\;  } s\to \infty$.
\end{itemize}
The slowly varying functions
$\sigma(s)=1+\ln s$, $s> 1$ (utility
function),
$\sigma(x)=S_0\,s/(S_0+s-1)$, and $S_0>1$ (sigmoidal function) are
examples
of $u$- and $s$-demands.
These functions describe various kinds of
diminishing returns.

Because the contribution from junior partners, $B(n,s)$, vanishes
as $s\downarrow s_0$, to ensure consistency in equations
(\ref{income}) and (\ref{stat}), equation (\ref{eqS0}) is
rewritten as follows:
\begin{align}
\label{minimal}
s_0-P(1)+s_0\, \sigma(s_0)\,R_0=0
\end{align}
and consider the equation to be the definition of $s_0$. Equations
(\ref{stat}) and (\ref{minimal}) represent the analytical setup of the model. The aim is
to describe the asymptotic behavior of $n(s)$ as $s\to \infty$ and
$s\to s_0$ and to prove the uniqueness of the solution. For
simplicity of notation, the price scale is chosen such that
$s_0=1$.

\subsection{Consequences of the continuity conjecture}
To study the problem on a
compact set,
the variable $s\in[1,\infty)$ is changed to
$x=1/s\in(0,1]$  in (\ref{stat}) and
(\ref{minimal}). According to
(\ref{minimal}), the integral in the
left-hand side of (\ref{stat}) vanishes
at the point $s=x=1$. Hence, the
right-hand side of (\ref{stat}) can be rewritten as the
following
product:
\begin{align}
\label{C}
s-P(s)+s\,\sigma(s)\,R_0=\sigma(s)\,(s-1)\,C(s^{-1})=
x^{-1}\sigma(x^{-1})(1-x)\,C(x),
\end{align}
where $\sigma(x^{-1})=x^{-1}$ for
$x=O(1)$.

As $s\to \infty $ and $x=1/s\to 0$, the
leading terms in the left- and
right-hand sides of (\ref{C}) are
$x^{-1}\sigma(x^{-1})C(0)$ and
$x^{-1}(1+\sigma(x^{-1})\,R_0)$, respectively.
Therefore, for
$C(0)=R_0+\sigma(\infty)^{-1}$,
we have
\begin{align}
\label{c0}
C(0)=\left\{
\begin{array}{ll}
R_0 &\text{for
$\ell$- and $u$-demands,}
\\
R_0+S_0^{-1} &\text{for
$s$-demand,}
\end{array}
\right.
\end{align}
because $\sigma(\infty)=\infty$ for
$\ell$ and $u$-demands, and $\sigma(\infty)=S_0$
for $s$-demand.
As $x\to 1$, we obtain
\begin{align}
\label{C1}
C(1)=1+2R_0-P'(1)\ge
1+2R_0>0,
\end{align}
where
$P'(1)\le0$ because the welfare $P(s)$
is normally a
decreasing (or finite) function.
Therefore, the left-hand side of
(\ref{C}) is an increasing function
that vanishes only at $s=1$.
For this reason, it is supposed that
$c_*>0$ exists such that
\begin{align*}
c(x)\defn {C(x)}/{C(0)}\ge c_* \;\; \forall
x\in(0,1].
\end{align*}
Suppose that $P(x^{-1})$ vanishes for
sufficiently small $x$, say for $x\in
(0, x_*]$, and $\sigma(s)=s$ for
$s\in[1,1/x_*]$ \footnote{A point
$x_*\in (0,1)$ will appear in some
of the assumptions below. In fact, these
additional assumptions can be imposed
at different points, but a fixed
point is used for simplicity of notations.}. It
is clear that for $x\in (0, x_*]$,
the assumption $\sigma(s)\le s$ implies
\begin{align}
\label{estC}
\frac{(1-x)C(x)}{R_0}=1+\frac{1}{\sigma(x^{-1})R_0}\ge
\left\{
\begin{array}{ll}
1+x/R_0 &\text{for
$\ell$- and $u$-demands,\;\;}
\\
1+1/S_0 R_0 &\text{for
$s$-demand.}
\end{array}
\right.
\end{align}
Let us search for a positive solution of
the balance equation (\ref{stat}) in
the factorized form
\begin{align}
\label{anzatz}
n(1/x)=x^b\,(1-x)^d\,\rho(x)/\sigma(x^{-1}),\quad
x\in(0,1], \quad b>1,\; d>-1.
\end{align}
As the normalization of $n$  is
not essential, $\rho(1)=1$ is used. Let us
derive an equation for $\rho$ and
calculate {\it a priori} values of $b$ and $d$. Set
\begin{align}
\label{Mn}
r(s)\defn R(s)/R_0,\quad
m_a\defn \int_1^\infty s^a\,
r(s)\,ds=1,\quad m_0=1
\end{align}
and suppose that  the probability density
$r(x)$ is a bounded continuous
function, $r(1)>0$.

\begin{theorem}
\label{contin} The relation $d={R(1)}/{C(1)}-1$
is necessary for the
continuity of a strictly positive
bounded function $\rho(x)$ at the point
$x=1$ for  all continuous demand functions such that $\sigma(1)=1$.
Continuity   at the point  $x=0$ implies
$b_\ell=3$,
$b_u=2$, and $b_s>2$
for $s$-demand; more precisely, $b_s$
satisfies the transcendental equation
\begin{align}
\label{b_s}
b_s:\quad
\int_1^{\infty}s^{b_s-2}r(s)\,ds=1+\frac{1}{S_0\,R_0}.
\end{align}
\end{theorem}

\proof According to equations (\ref{C}) and (\ref{anzatz}),
rewrite equation (\ref{stat}) in  the variables $x=1/s$, $y=1/z$:
\begin{align*}
n(s)\,\bigl(s-P(s)+s\,\sigma(s)\,
R_0\bigr)=&R_0\,x^{b-1}(1-x)^{d+1}c(x)\,\rho(x)
\\
&=
\int_{1}^s\sigma(y)
\,R(s/y)\,n(y)\,dy=\int_{1/s}^1
\sigma(z^{-1})\,z^{-2}R(sz)\,n(z^{-1})\,dz
\\
&=
\int_{x}^1\,z^{b-2}\,
(1-z)^d \,R(z/x)\,\rho(z)\,dz.
\end{align*}
This integral equation for
$\rho(x)$
can be represented in two  equivalent forms which are indepemdent
of $\sigma$:
\begin{align}
\label{eqrho}
\rho(x)=&\frac{1}{x(1-x)C(x)}\int_x^1
\biggl(\frac{z}{x}\biggr)^{b-2}\,
\biggl(\frac{1-z}{1-x}\biggr)^d
R(z/x)\,\rho(z)\,dz
\\
&
\label{eqrho2}
=\frac{1}{(1-x)C(x)}\int_1^{1/x}
y^{b-2}\,
\biggl(\frac{1-xy}{1-x}\biggr)^d
R(y)\,\rho(xy)\,dy.
\end{align}
If $d>-1$ and $x>0$ is fixed,
the measure
$\mu(x,dz)=(d+1)(1-z)^d/(1-x)^{d+1}\,dz$
is a  probability
measure on $(x,1]$. If $\rho(x)$ is a
bounded measurable function that is
continuous at the point $x=1$,  so is
the product
$$
p(x,z)=({z}/{x})^{b-2}
\,R(z/x)\,\rho(z),
\quad z\in(x,1)
$$
and
$\overline{\lim}\,p(x,z)=\underline{\lim}\,p(x,z)=
R(1)\,\rho(1)$ as $x\to 1$, $z\in (x,1]$.
Denote with $\e_x$ the mathematical expectation with respect to
the probability  measure $\mu(x,dz)$. Hence,
$$
\overline{\lim}\,\e_x\,
p(x,\cdot)=\underline{\lim}\,\e_x\,
p(x,\cdot)= R(1)\,\rho(1).
$$
Multiplying (\ref{eqrho}) by
$d+1$,  we obtain
the condition that is necessary for  the continuity of $\rho$:
\begin{align}
(d+1)\,\rho(1)=\lim_{x\to
1} \e_x\, p(x,\cdot)/({x
\,C(x)})=R(1)\,\rho(1)/C(1).
\end{align}
If $\rho(1)\neq 0$, this equation
yields (\ref{d}):
$d=R(1)/C(1)-1=r(1)/c(1)-1> -1$.

As $x\to0$,  the left-hand side of  (\ref{eqrho2})
clearly converges to $\rho(0)$;
the Lebesgue theorem on dominated convergence
implies that
the right-hand
side of equation (\ref{eqrho2}) converges to
$$
\frac{\rho(0)\,R_0}{C(0)}\,m_{b-2}=\rho(0),\quad
m_\alpha\defn\int_1^{\infty}s^{\alpha}
r(s)\,ds.
$$
As $C(0)/R_0=1$ for $\ell$- and
$u$-demands,
and   $C(0)/R_0=1+1/S_0R_0$ for $s$-demand,  the
necessary continuity condition at $x=0$ is
the equation with respect to $b$
\begin{align}
\label{b}
m_{b-2}=
\left\{
\begin{array}{ll}
1 &\text{for
$\ell$- and $u$-demands,\;\;}
\\
1+1/R_0S_0 &\text{for
$s$-demand}
\end{array}
\right.
\end{align}
provided $\rho(0)\neq0$
(see (\ref{c0}) and (\ref{Mn})). The moment $m_a$ is an increasing
function of $a$ because $r$ is a regular probability density on
$[1,\infty)$; hence, equations (\ref{b_s}) or
(\ref{b}), in each case, have a unique solution.
\hfill $\square$
\medskip

{\sl Remark 1.}
As $\sigma((xy)^{-1})/\sigma(x^{-1})\to 1$
as $x\to0$,
this result remains valid
if $\rho(x)=g>0$.
This will be the case for $u$-demand. It is also clear that the
Pareto exponent for $n(s)O(\sigma( x^{-1})^{-g})$ is the same as that for
$n(s)$.

{\sl Remark 2.}
If $\delta\defn 1/S_0 R_0$  is a small constant and the measure
$r(s)\,ds$ has at least one moment, the
equation $m_{\alpha}=1+\delta$ can be
approximately solved:
$$
\alpha=\frac{\delta}{\int_1^\infty
r(s)\,\ln s\, ds}+O(\delta^2)>0,\quad \delta= 1/S_0
R_0,
$$
i.e., $a=1+\alpha$, similar to the Yule
case\footnote{The constant $R_0$
characterizes an
intensity of pairwise interactions.
In this context,
$1/R_0$ is the average time between
interactions. On the other hand, $S_0$ can be interpreted as
the volume of a  ``market''.  The
difference $a-1$ between the Pareto exponent
$a$ ($a>1$) and the exponent of the
``heaviest'' tail is termed the {\it horizon}
of the distribution.
In these
terms, the Pareto exponent
that corresponds to the sigmoidal demand satisfies
the following law: ``{\it Horizon
$\times$ market volume $\thicksim$
expectation-of-transaction time}.''}.

Thus, it is demonstrated that the continuity and
strict positivity of $n(s)$ implies the
expansion:
$$
u(x^{-1})=x^{b}(1-x)^{d}\rho(x)/\sigma(x^{-1}),\quad
b-1=a_{\rm net}=\left\{
\begin{array}{ll}
1 &\text{for
$\ell$- and $u$-demands,}
\\
1 +\alpha &\text{for
$s$-demand.}
\end{array}
\right.
$$
In addition, for all kinds of demands, equation (\ref{eqrho})
can finally be rewritten as follows:
\begin{align}
\label{eqrhoS}
\rho(x)=\frac{1}{x(1-x)C(x)}\int_x^1
\biggl(\frac{z}{x}\biggr)^{\alpha}
\biggl(\frac{1-z}{1-x}\biggr)^d
R(z/x)\,\rho(z)\,dz,
\end{align}
where $\alpha=0$ for $\ell$- and
$u$-demands, and $C$ is the only function
that is dependent on $\sigma$ (see equation (\ref{C})).

\subsection{Pareto exponents for gross income}
The relationship $g=s+\sigma(s)\,R_0$
between gross and net income allows
one to calculate the Pareto exponents
for the gross income distribution.
First consider the  $\ell$-demand. Denote
the unique  positive root of the
equation $g=s+s^2\,R_0$ with $s(g)$, where $g$ is gross income. The density
distribution of  $g$ then equals
$N_\ell(g)$:
$$
N_\ell(g)=n_\ell(s(g))\,s'(g),\quad
s(g)=1/\sqrt{1+4R_0\,g}=O(g^{1/2}),
$$
where $n_\ell(s)=O(s^{-3})$,
$n_\ell(s(g))=O(g^{-3/2})$, and $s'(g)=
O(g^{-1/2})$. Finally, we obtain
$N_\ell(g)=O(g^{-2})$, i.e., the Pareto
exponent  of  the gross income  $a_{\rm
gross}=a_{\rm net}-1=1$.

Now consider  the $u$-demand. Because
$\sigma(s)$ is a slowly varying
and increasing function,  there exists a
finite constant $\sigma_\varepsilon$
such that $\sigma(s)\le
\sigma_\varepsilon s^\varepsilon$ for
any $\varepsilon>0$. Now denote
the unique positive root of the
equation $g=s\,(1+\sigma(s)\,R_0)$ with $s(g)$, and
let $N_u(g)$ be the probability
density of $g$. The monotonicity
of $\sigma(s)$ implies the following inequalities:
\begin{align}
\label{ineq}
g\ge s(g)=\frac{g}{1+\sigma(s(g))\,R_0}\ge \frac{g}{1+\sigma(g)\,R_0}\ge
\frac{g}{1+\sigma_\varepsilon
g^\varepsilon\,R_0}.
\end{align}
By definition of the cumulative probability, we have
\begin{align*}
\int_y^\infty\,N_u(g)\,dg=\int_y^\infty
n_u(s(g))\,s'(g)\,dg
=\int_{s(y)}^\infty\,n_u(g)\,dg.
\end{align*}
As $n_u(g)=O(g^{-2})$,
the above inequality (\ref{ineq}) yields the
two-sided
estimates for the cumulative probability $\p_u\{g\ge y\}$:
\begin{align*}
O(y^{-1})=&\int_{y}^\infty\,n_u(g)\,dg
\le \int_{s(y)}^\infty\,n_u(g)\,dg=\int_y^\infty\,N_u(g)\,dg
=Z^{-1}\,\p_u\{g\ge y\}
\\
&\le
\int_{y/(1+\sigma_\varepsilon
y^\varepsilon\,R_0)}^\infty\,n_u(g)\,dg=
y^{-1}O(1+\sigma_\varepsilon
y^\varepsilon\,R_0),
\end{align*}
where $Z$ is  the  normalizing constant and
$\varepsilon$ is arbitrarily small. Hence,
passing to the limit as $y\to\infty$,
$\varepsilon\to0$, we have
$$
a_{\rm gross}=-\overline{\lim}\,\ln\bigl(\p\{g\ge y\}\bigr)/\ln
y=-\underline{\lim}\,\ln\bigl(\p\{g\ge y\}\bigr)/\ln
y=1=a_{\rm net}.
$$

For $s$-demand, gross
and net income have the same order:
$g(s)=s\,(1+\sigma(s)\,R_0)=s+O(s)$,
as $\sigma(s)$ is  bounded by
$S_0$. Hence,
\begin{align}
\label{alpha}
a_{\rm net}=a_{\rm
gross}=1+\alpha, \qquad \alpha:\;1+\frac{1}{S_0R_0}=\int_1^\infty
s^\alpha\,r(s)\,ds.
\end{align}

The above estimates of $b$ and $d$
define meaningful asymptotical
properties of $n(s)$ provided that $\rho(x)$
is bounded from above and below. The proof of
this fact  is a nontrivial mathematical
problem that will be solved via a
probabilistic representation of the
solution of equations
(\ref{stat}) and (\ref{minimal}); however, first
the results of an empirical
study of this problem are considered.

\section{Numerical simulation of $\rho$}
With respect to the analytical proof of the boundedness and
the strict positivity of  $\rho(x)$
etc.,  a computer simulation
is an easier approach that
convinced us that equation
(\ref{eqrhoS}) is
consistent.
Although the operator $A:\,C[0,1]\to C[0,1]$
\begin{align}
\label{OpA}
(A\rho)(x)=\frac{1}{x(1-x)c(x)}\int_x^1
\biggl(\frac{z}{x}\biggr)^\alpha\,
\biggl(\frac{1-z}{1-x}\biggr)^{r(1)/c(1)-1}
r(z/x)\,\rho(z)\,dz
\end{align}
is not a contraction in the uniform
topology\footnote{The estimate  given in Remark 3 implies that
$||A||_{C}\le F_*\Pi_1\Pi_2<\infty$, and  the numerical
spectral estimates show that
$||A||_{C}\ge 1$.}, the iteration procedure
$
\rho_n=A\rho_{n-1}
$
converges to a smooth and strictly positive
solution $\rho(x)$ for all tested probability distributions $r(s)$,
initial approximations $\rho_0$, and
values of $d\in(-1,\infty)$.
Figures 2-4  show  the numerical
solutions $\rho(x)$ and $n(s)$ of equation
(\ref{eqrhoS}) for the pair-interaction
$R(s)=\lambda
e^{-\lambda(s-1)}$ and ``earning''
$P(s)=2/(1+(s-1)^3)$;  black, red and blue curves
correspond to  $d=0$,  $d>0$, and
$d<0$, respectively.

%\psdraft
\begin{figure}

\begin{minipage}{0.5\linewidth}
%\centering\epsfig{figure= S0x0p0,width=\linewidth}
\centering\epsfig{figure= 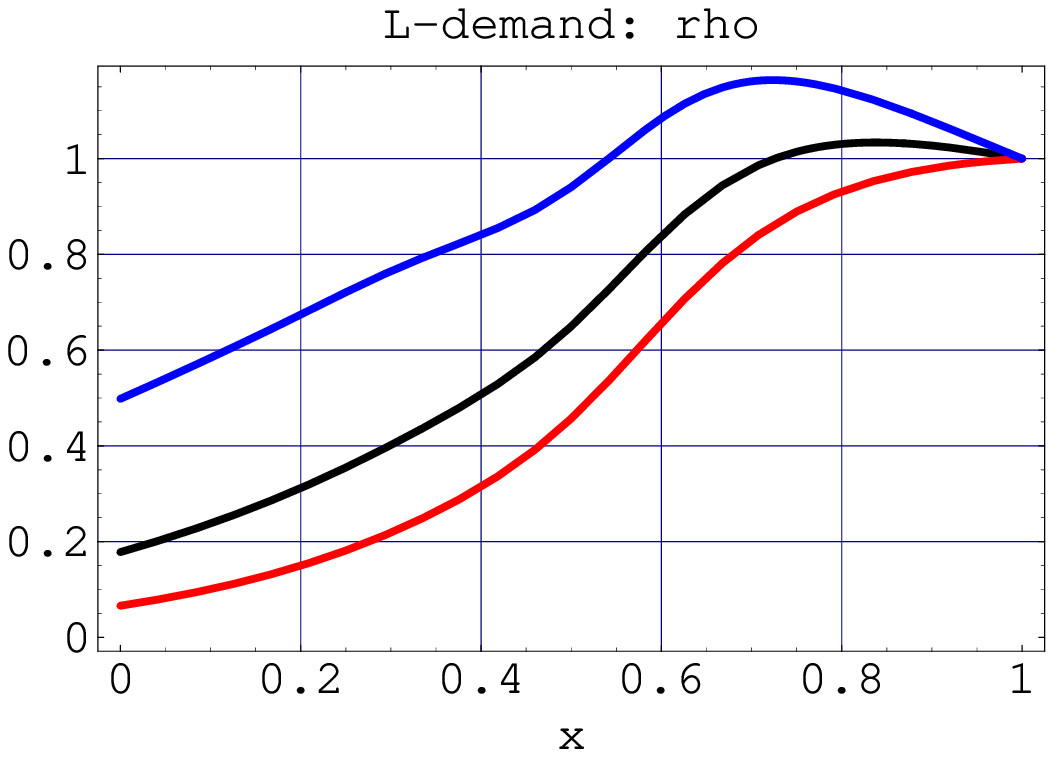,width=\linewidth}
\end{minipage}\hfill
\begin{minipage}{0.5\linewidth}
%\centering\epsfig{figure= S0x0p0around0.EPS,width=\linewidth}
\centering\epsfig{figure= 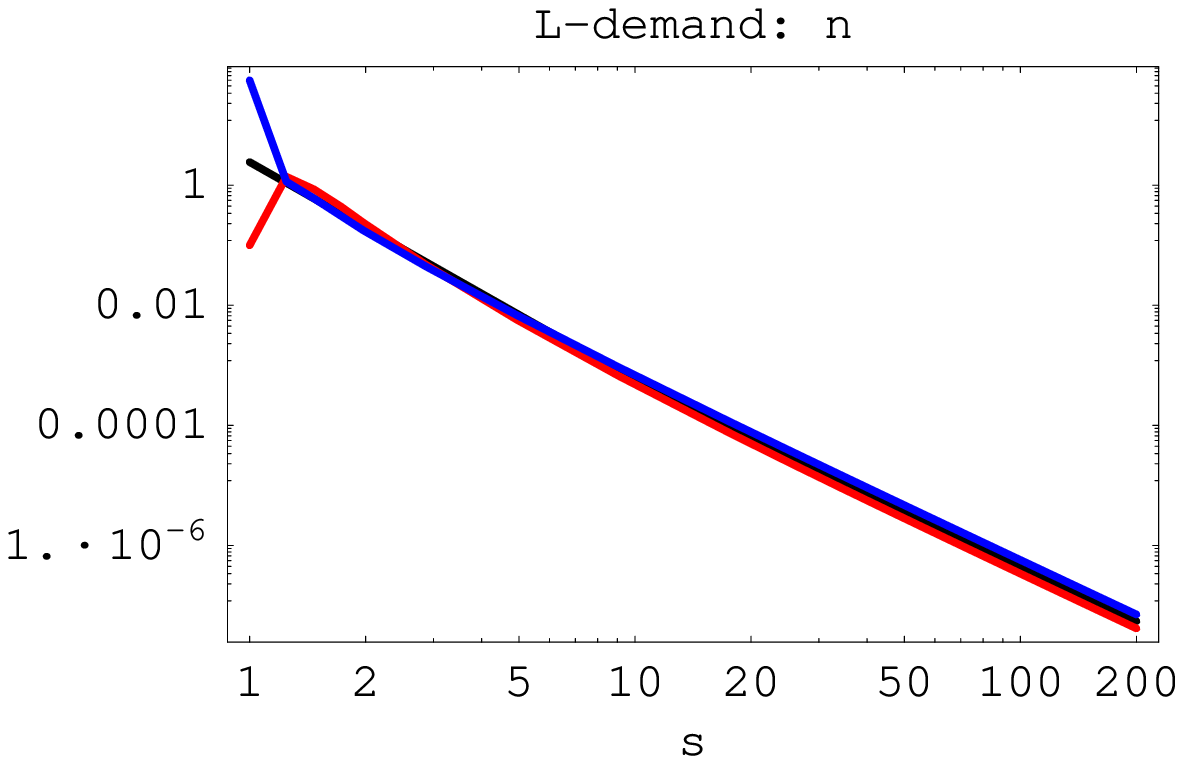,width=\linewidth}
\end{minipage}\hfill
\caption{\small For $\ell$-demand, we set
$\sigma(s)= s$ and $n(s)=n_0\,s^{-3}(1-s)^d\rho(1/s)$;
more precisely,
$\lambda=3$ for $d=0$ (black curve),
$\lambda=4$ for
$d=1/3$ (red curve), $\lambda=2$
for $d=-1/2$ (blue curve).}
\vskip3mm

\begin{minipage}{0.5\linewidth}
%\centering\epsfig{figure= S0x0p0,width=\linewidth}
\centering\epsfig{figure= 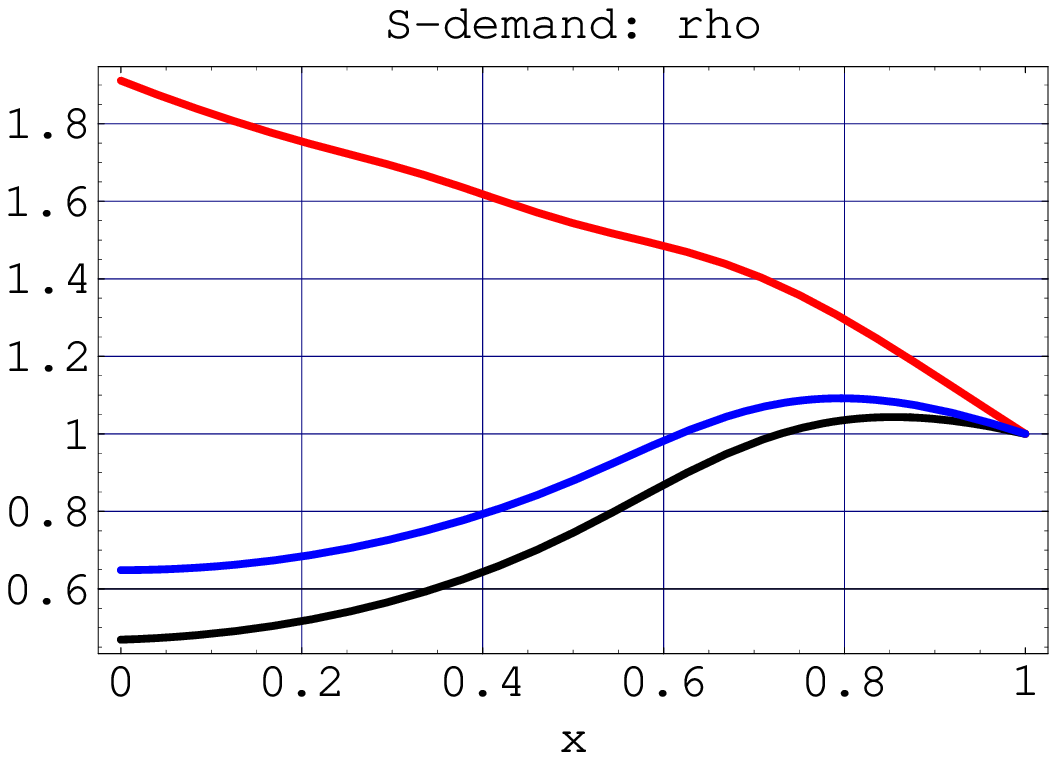,width=\linewidth}
\end{minipage}\hfill
\begin{minipage}{0.5\linewidth}
%\centering\epsfig{figure= S0x0p0around0.EPS,width=\linewidth}
\centering\epsfig{figure= 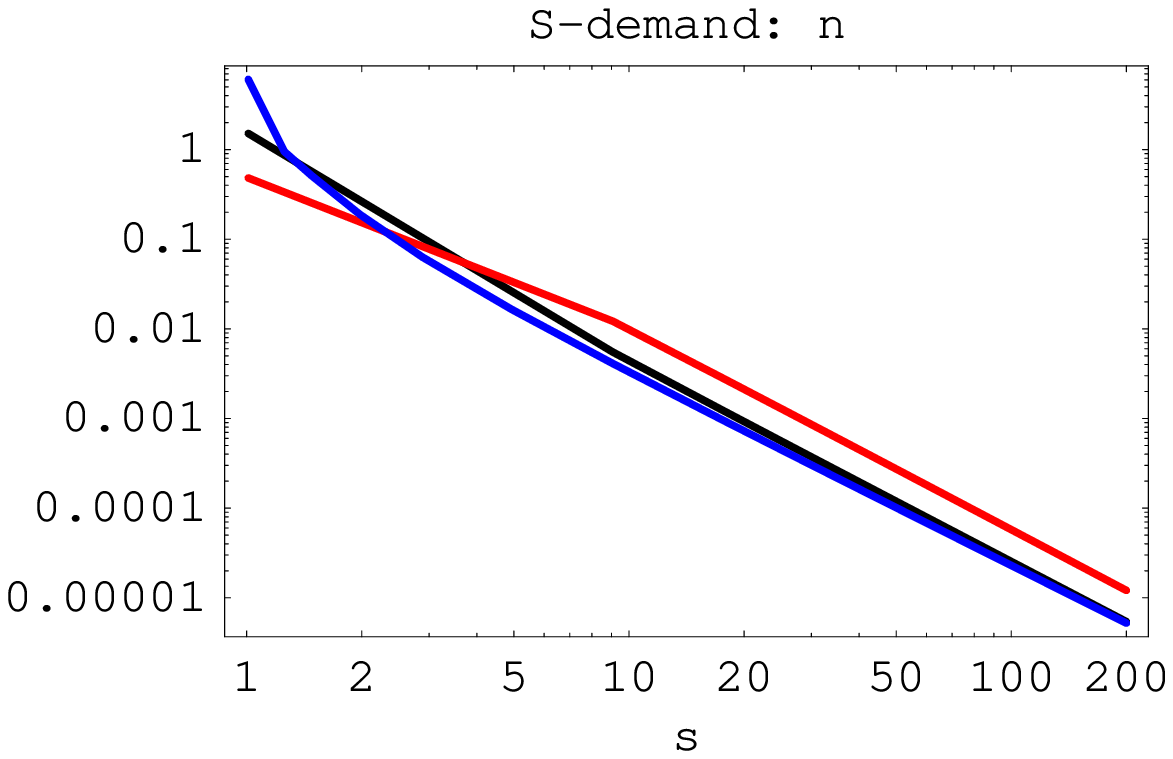,width=\linewidth}
\end{minipage}\hfill
\caption{
\small
In order to illustrate $\rho$ and $n$ for $\sigma$-demand, we
set
$\sigma(s)=20/(1+19s^{-2})$.
Solving equation (\ref{b_s}), we find
$b = 2.22741$ for $d=0$ and $\lambda=3.9$,
$b = 2.25672$ for
$d=1/39$ and $\lambda=4$, $b=2.13304$
for $d=-19/39$ and $\lambda=2$, so that
$n(s)=n_0\,s^{-b}(1-s)^d\rho(1/s)/\sigma(s)$.}

\vskip3mm
\begin{minipage}{0.5\linewidth}
%\centering\epsfig{figure= S0x0p0,width=\linewidth}
\centering\epsfig{figure= 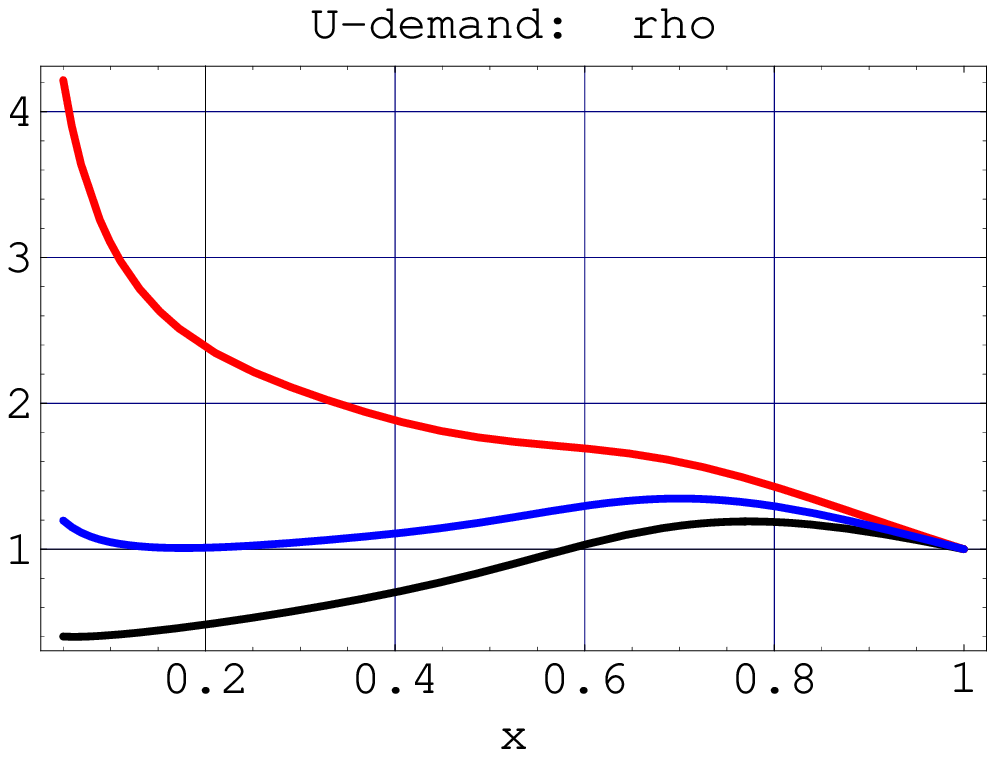,width=\linewidth}
\end{minipage}\hfill
\begin{minipage}{0.5\linewidth}
%\centering\epsfig{figure= S0x0p0around0.EPS,width=\linewidth}
\centering\epsfig{figure= 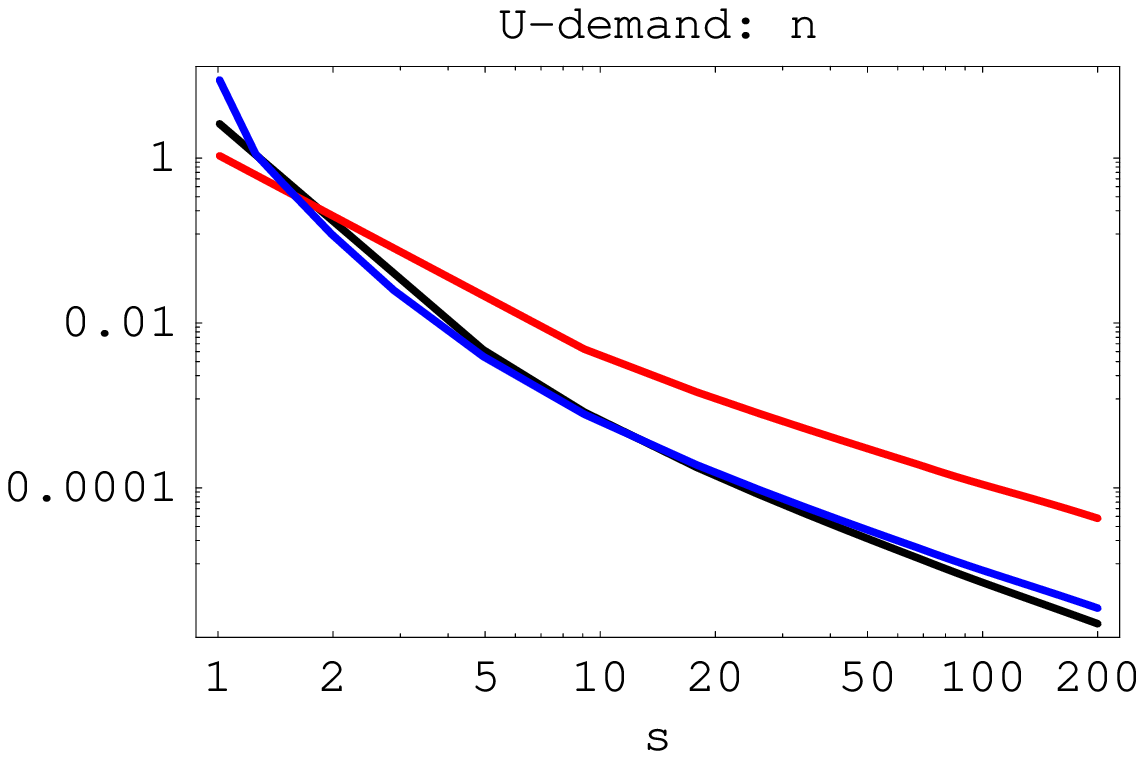,width=\linewidth}
\end{minipage}\hfill
\caption{\small For $u$-demand, we have
$\sigma(s)=1+\log s$ and
$n(s)=n_0\,s^{-2}(1-s)^d\rho(1/s)/\sigma(s)$.
Simple computation implies
$\lambda=3$ for $d=0$,
$\lambda=4$ for
$d=0.3$, and $\lambda=2$
for $d=-1/3$.}
\end{figure}

It is worthwhile noting that the solution for equation
(\ref{eqrhoS})  can be
calculated on any interval $(y,1]\subset (0,1]$. This
means that problem (\ref{stat})  remains
well-posed on any finite interval $[s_{\rm min},s_{\rm
max}]$ and $n(s_{\rm max})=O(s_{\rm
max}^{-b})$.

It was observed that the Fourier or
Chebyshev polynomial approximation of
the solution $\rho(x)$ does not provide
a fast convergence. By differentiating
equation (\ref{eqrhoS}), one can derive
a pair of homogeneous boundary
conditions of the form
$\rho'(x)=\gamma_x\rho(x)$, $x=0,1$.
Note that the standard orthogonal
polynomials do not fulfill  the
boundary conditions of this type; this
is a possible explanation of
the inefficiency of the polynomial
approximation. Nevertheless, the finite
dimensional approximation of operator
(\ref{OpA}) in the Fourier basis shows
that the spectrum $\{\lambda_k\}$ of
$I-A_n$ is distributed more or less
uniformly over a smooth ``semicircle''
curve  ($\re \lambda_k\in[-1,0]$, $\im
\lambda_k\in[0,1/2]$)  (see Fig. 5).

\begin{figure}
\begin{minipage}{0.5\linewidth}
%\centering\epsfig{figure= S0x0p0,width=\linewidth}
\centering\epsfig{figure= 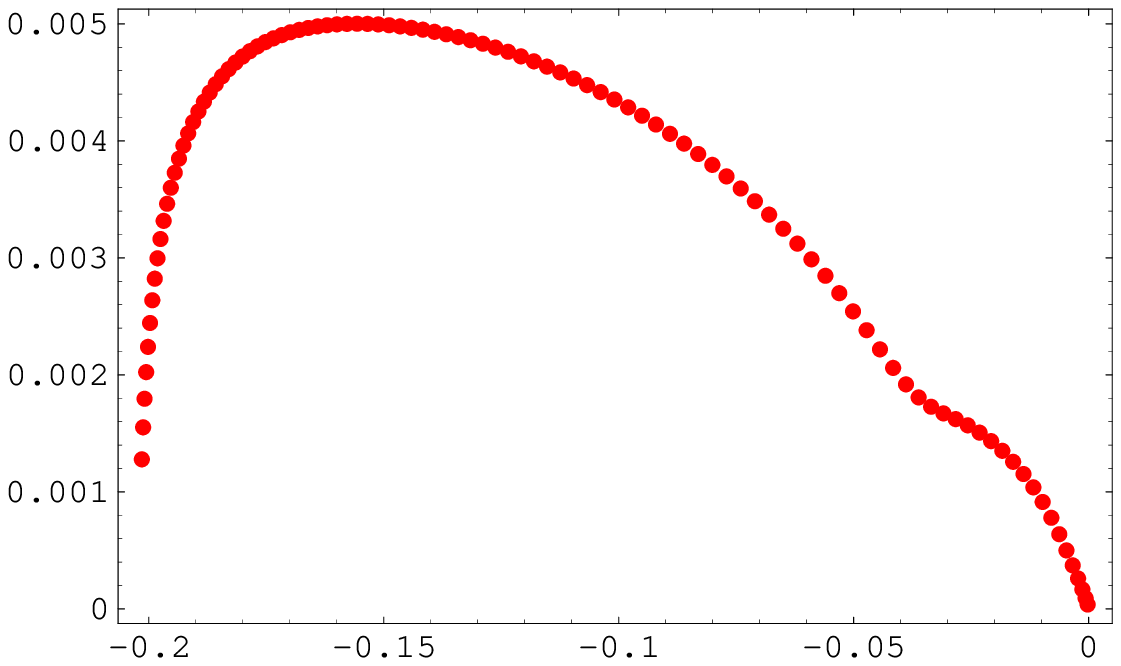,width=\linewidth}
\end{minipage}\hfill
\begin{minipage}{0.5\linewidth}
%\centering\epsfig{figure= S0x0p0around0.EPS,width=\linewidth}
\centering\epsfig{figure= 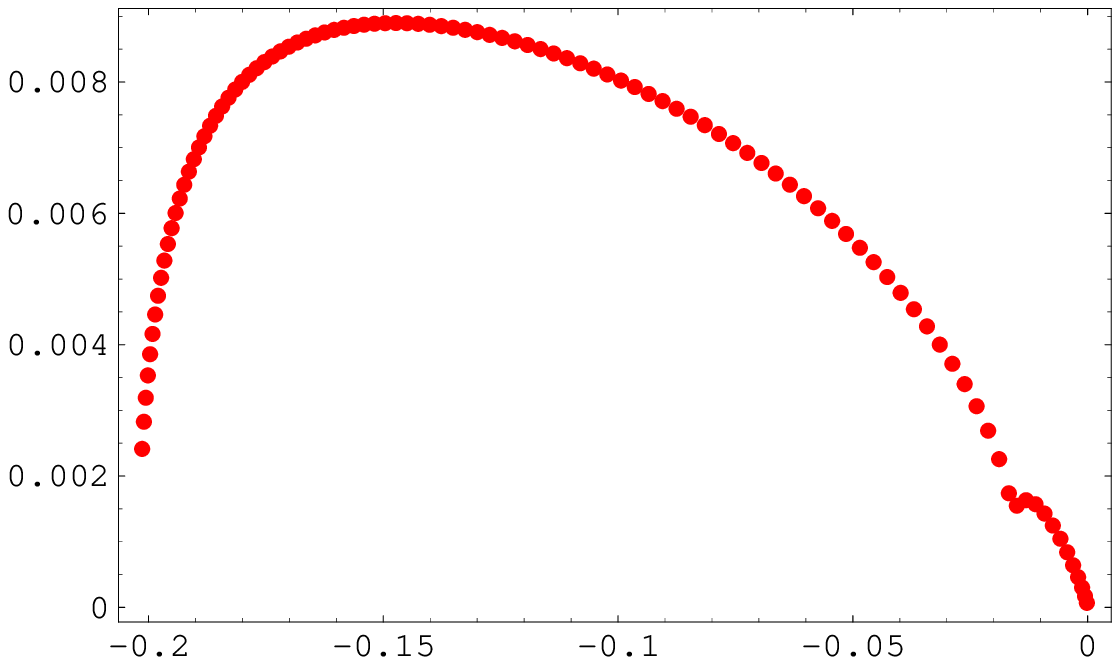,width=\linewidth}
\end{minipage}\hfill
\caption{\small The spectrums of the operator $I-A_n$ in
$\ell$-cases for $d=0$ and $d=-0.5$,
in the complex plain for $n=100$. There is no reason to expect the
existence of a spectral gap; it is most likely that the limit spectrum is
continuous.}
\end{figure}

The spectrum has no spectral gap near
the minimal eigenvalue
$\lambda_0\approx 0$, however, in
a finite dimensional approximation at least, the
minimal eigenvalue is simple. Thus, the
spectral approximation provides no helpful
hints as to the proof of the uniqueness
theorem based on the existence of a
spectral gap. In contrast, Fig. 4
shows that all eigenvalues $\lambda$,
such that $|1-\lambda|>1$, are complex;
the corresponding eigenfunctions
must also be complex. This observation
illustrates and extends the statement
which will be
proved in
Theorem \ref{Uniq} below: a positive
bounded solution of equation
(\ref{stat}) is unique, while all other
solutions are complex. This
statement is a kind of
Perron-Frobenius theorem.

Now we must pay attention to the
interesting mathematical background to the
problem (\ref{stat})-(\ref{minimal}).
The consideration below is
self-consistent and does not refer the
reader to mathematical
literature.

\section{Markov chains
and {\it a priori}  bounds for
$\rho(x)$}
\subsection{Transition probabilities near $0$ and $1$}

Consider the Markov chain $\{x_n\}$ on
$(0,1]$ with  transition probability
$P(x,A)$ from $x$ to $A\subset [x,1]$
given by a pair of measures describing
the process in the neighborhoods of the
points $0$ and $1$:
\begin{align}
\label{MC}
P(x,A)=\left\{
\begin{array}{cl}
\frac{1}{x\,Z_\alpha(x)}
\int_A \,(y/x)^\alpha \,r({y}/{x})\,dy,
& \mbox{ if }
x\in(0, x_*],
\\
 {}
 \\
 -\int_A\,
d \left(\frac{1-y}{1-x}\right)^{d+1},
& \mbox{ if } x \in (x_*,1],
\end{array}
\right.
\end{align}
where $\alpha=0$ for $\ell$- and
$u$-demands; for $s$-demand, $\alpha$
satisfies equation (\ref{alpha}); and $Z_\alpha$
is a
local
normalizing factor:
\begin{align}
\label{Z}
Z_\alpha(x)\defn
\int_{[1,1/x]}\,s^\alpha \,r(s)\,ds \le
1.
\end{align}
Note that $Z_0(x)\le 1$ and $Z_\alpha(x)\le 1+1/S_0R_0$
(see (\ref{estC})), but  $Z_\alpha(x)/(1-x)c(x)\le
1$ in any case on
$(0,x_*]$.

The associated
function $F(x,y)$ on the set $x,y: \,1> y\ge x$
$$
\label{F}
F(x,y)=\left\{
\begin{array}{cl}
\frac{Z_\alpha(x)}{(1-x)\,c(x)}
\left(\frac{1-y}{1-x}\right)^{d},
& \mbox{ if }
x\in [0,x_*],
\\
 {}
 \\
 \frac{1}{x}\,
\frac{c(1)}{c(x)}\,\frac{r(y/x)}{r(1)}\,
\left(\frac{y}{x}\right)^\alpha,
& \mbox{ if } x \in (x_*,1],
\end{array}
\right.
$$
is defined such that equation
(\ref{eqrhoS}) is
equivalent to the integral equation
\begin{align}
\label{expect}
\rho(x)=\e_xF(x,x_1)\rho(x_1),\quad
x\in(0,1].
\end{align}
As the mathematical expectation $\e_x$,
$x\in(0,1]$,
is the well-defined integral operator,
any measurable bounded solution $\rho(x)$ is
absolutely continuous; moreover,
equation
(\ref{expect}) can be rewritten as
$$
\rho(x)=\e_x\prod_{k=1}^n
F(x_{k-1},x_k)\;\rho(x_n).
$$
According to Lemma 1 from the Appendix,
the Markov chain exits any interval $(x,y)\subset(0,1]$
with probability one. Hence, if the
solution $\rho$ is bounded, it is continuous and
$\rho(x_n)\to \rho(1)=1$ as
$n\to\infty$ with probability one;
therefore,
\begin{align}
\label{repr}
\rho(x)=\e_x\prod_{n=1}^\infty
F(x_{n-1},x_{n}),\quad x_0=x,
\end{align}
where $\{x_n\}_1^{\infty}$ is  the
monotone discrete-time Markov chain
to be defined by (\ref{MC}).

Consider a priori bounds for
mathematical expectation (\ref{repr}).
For simplicity, it is supposed that the
probability density $r(s)$ is uniformly
bounded. Then $x\le y$ on the domain of
$F(x,y)$.  It is clear that $F(x,y)\to
1$, if either $y\to 0$ or $x\to1$. Set
$$
F_*=\max_{0\le x\le x_*\le
y\le 1}F(x,y)
$$
and consider the integer-valued exit
time $n_*=n(x,x_*)\ge 0$ from the
segment
$(0,x_*]$ for the Markov chain
$\{x_n\}_{n=0}^{N_*}$ staring at
$x_0=x$:
$x_{n_*-1}\le x_*<x_{n_*}$
. We have
\begin{align*}
\prod_{n=1}^\infty
F(x_{n-1},x_n)=& F(x_{n_*-1},x_{n_*})\prod_{n=1}^{n_*-1}
F(x_{n-1},x_n)
\prod_{m=1}^{\infty}
F(y_{m-1},y_m)
\\
&\le F_*
\prod_{n=1}^{n_*-1}
F(x_{n-1},x_n)\,
\prod_{m=1}^{\infty}
F(y_{m-1},y_m),
\end{align*}
where  $\{y_m\}$ is the Markov chain
starting
at the point $y=y_0=x_{n_*}\in[x_*,1]$. Thus,
$$
\e_x\prod_{n=1}^\infty
F(x_{n-1},x_n)\le
F_*\,\e_x\prod_{n=1}^{n_*-1}
F(x_{n-1},x_n)\,\sup_{y\ge
x_*}\e_y\prod_{m=1}^\infty
F(y_{m-1},y_m).
$$
Therefore, it suffices to prove  the
uniform boundedness of
$\Pi_1(x)=\e_x\prod_{n} F(x_{n-1},x_n)$
for $x_k\in(0,x_*]$ and $\Pi_2(x)=\e_x\prod_{m}
F(y_{n-1},y_n)$ for $y_m\in [x_*,1)$.

\subsection{Upper bound and the uniqueness theorem}
Consider the  uniform upper
bound for $\Pi_1$.
According to (\ref{estC}) and (\ref{Z}),
\begin{align*}
{Z_\alpha(x)}/{(1-x)c(x)}&\le1  \quad \forall x\le x_*.
\end{align*}
The Markov chain $\{x_n\}$ increases
monotonically and exits any interval
$(x,y)\subset(0,x_*]$ with  probability
one (see Appendix, Lemma 1). Thus, we have
$$
\Pi_1(x)=\e_x\prod_{n=1}^{n_*-1}
F(x_{n-1},x_n)\le
\biggl(\frac{1-x_{n_*-1}}{1-x}\biggr)^{d}
\le
\left\{
\begin{array}{cl}
1,
& \mbox{ if } d\ge0,
 \\
(1-x_*)^{d},
& \mbox{ if }
d< 0
\end{array}
\right\}
<\infty.
$$

Consider the  upper bound for $\Pi_2$.
Note that $F(1,1)=1$,
and suppose that
\begin{align}
\label{F^*}
F^*\defn &\max_{x\ge
x_*,\,y\in[x,1]}\frac{1}{1-x}\,\bigl|1
-F(x,y)\bigr|
\nonumber
\\
&=
\max_{x\ge
x_*,\,y\in[x,1]}\frac{1}{1-x}\,\biggl|1
-\frac{c(1)}{x\,c(x)}\frac{r(y/x)}{r(1)}
\left(\frac{y}{x}\right)^\alpha\biggr|<\infty.
\end{align}
For example, if $R(s)$ is a monotone
function and $\alpha=0$, then $r(y/x)\le r(1)$, and
\begin{align}
\label{F*}
F^*\le \frac{1}{x_*\,c_*}\biggl(c(1)+\max_{x\in[x_*,1)}
\frac{|c(1)-c(x)|}{1-x}\biggr),
\end{align}
where $c_*$ is $\min c(x)$.
In any event,  $F(y_{n-1},y_{n})\le 1+(1-y_{n-1})\,F^*$, and the estimate
$$
\Pi_2(y)=\e_y\prod_{m=1}^\infty
F(y_{m-1},y_m)\le  \e_y \exp
\{F^*\,\sum_{m}(1-y_{m})\},\quad y_0=y
$$
can be rewritten in terms of  the
independent variables $\xi_n\defn (1-y_n)/(1-y_{n-1})\in(0,1]$ as
the sum of a stochastic geometrical progression:
$$
\sum_{m=0}^{\infty}(1-y_{m})=(1-y)\biggl(1+
\sum_{m=1}^{\infty}\prod_{k=1}^m\xi_n
\biggr),\quad \xi_m\in(0,1].
$$

\begin{theorem}
\label{UpperBound}
Let  $F^*$ be defined by {\rm(\ref{F^*})}.
Then,
$\rho(x)$ is a bounded
continuous function if
\begin{align}
\label{est}
d+1<F^*\,\frac{2-e^{F^*}}{e^{F^*}-1-F^*},
\quad F^*<\ln 2.
\end{align}
\end{theorem}

\begin{proof}
As proved above, $\sup_{x\in(0,x_*]}\Pi_1(x)<\infty$
and $F_*<\infty$
under rather general assumptions.
To prove that $\sup_{y\in(x_*,1)}\Pi_2(y)<\infty$,
it suffices to ensure
the finiteness of the exponential
moment
\begin{align}
\label{geom}
\Pi_2(y)=\e_y\,e^{(1-y)F^*\,(1+G)}<\infty, \quad
G=\sum_{m=1}^{\infty}\prod_{k=1}^m\eta_n.
\end{align}
According to (\ref{MC}), the factors
$\eta_m\in[0,1]$ of  the random geometrical
progression (\ref{geom}) are  independent random
variables with the following probability distribution:
\begin{align}
\label{xi}
P(\eta_m\in
B)=\gamma\int_B\eta^{d}\,d\eta,
\quad \gamma=d+1> 0.
\end{align}
Note that $\mu_n=\gamma/(\gamma+n)$ for  random variable
(\ref{xi}), and that according to
(\ref{suff1}) from Lemma 2 (see Appendix),  the exponential moment
of the geometric progression
(\ref{geom}) is finite if
$$
\sup_{k\ge 1}\sum_{n=1}^\infty
\frac{\lambda^n}{n!}\,\frac{1}{1-\mu_{n+k}}=
\sum_{n=1}^\infty
\frac{\lambda^n}{n!}\,\biggl(1+\frac{\gamma}{n+1}\biggr)\le
e^{\lambda}-1+\frac{1}{\lambda\gamma}\bigl(e^{\lambda}-1-\lambda\bigr)
<1
$$
for $\lambda=F^*$. This
inequality can be
resolved with respect to $\gamma$:
$\gamma<\lambda\,(2-e^\lambda)/(e^\lambda-1-\lambda)$.

Integral equation
(\ref{expect}) implies
the absolute continuity of a bounded
solution $\rho$.
\hfill $\square$
\end{proof}

\medskip
{\sl Remark 3}
Condition (\ref{est}) is the most
restrictive assumption.
For example, as $R(s)$ is a monotone function
in cases  (i--iii) considered in Section 3,
estimate (\ref{F*}) can be used
to evaluate $F_*=20$.
The convergence still holds
despite the fact that the  sufficient condition
$F_*\le \ln 2\approx 0.693$  (see (\ref{est})) is
crudely violated, i.e., the
sufficient condition is rather far from
an unknown necessary criteria.
It appears that assumption (\ref{est}) can be relaxed after a
more detailed study of exponential
moments.

\medskip
A remarkable fact to arise from this analysis is that a
solution, bounded from above,
is unique,
independently of the lower bound.

\begin{theorem}
\label{Uniq} The class of strictly
positive bounded continuous solutions
of equation {\rm(\ref{eqrho})} consists
of the unique function.
\end{theorem}

\proof Suppose that there exists two
nonnegative bounded continuous
solutions, $\rho_1(x)$ and $\rho_2(x)$,
that satisfy the condition $\rho(1)=1$.
The continuity assumption implies that
 $d_1=d_2$ and $b_1=b_2$ for
$n_1(s)$ and $n_2(s)$, respectively.

Consider the bounded and uniformly continuous
function $\Delta(x)=|\rho_1(x)-\rho_2(x)|$ on $[0,1]$.
It is  clear that
$$
\max_{x\in
[x_*,y]}P(x,y)=-\min_{x\in
[x_*,y]}\int_x^{y}
\,d\biggl(\frac{1-y}{1-x}\biggr)^{d+1}
=1-\biggl(\frac{1-y}{1-x_*}\biggr)^{d+1}<1
$$
majorizes the probability for the
Markov chain to remain in  the segment
$[x,y]$, $y<1$, $x\in[x_*,y]$ after
one step, and hence the  chain exits
any segment $[x_*,y]$ with a probability
of one (see also Lemma 1 in Appendix). In contrast, $\e_x\prod_{1}^{\infty}
F(x_{n-1},x_n)\le \Pi_2(x)<\infty$,
where $\Pi_2(x)$ is uniformly bounded
in $[x_*,1)$ (see Theorem 2). Equation
(\ref{expect}) implies that
$$
\Delta(x)\le \e_x\prod_{n=1}^{n_y}
F(x_{n-1},x_n)\max_{y\in[x_{n_y},1]}\Delta(y)\le
\Pi_2(x)\,\max_{\xi\in[y,1]}\Delta(\xi),
$$
where $n_y$ is the exit time from $[x_*,y]$:
$x_{n_y-1}<y\le x_{n_y}$.
As $y$ can
be arbitrarily chosen to be close to $1$
and
$$
\lim_{y\to1}
\max_{\xi\in[y,1]}\Delta(\xi)=0,
$$
(because of the uniform continuity)
the previous estimate implies that
$\Delta(x)=0$ for all $x\ge x_*$.
The uniform boundedness of $\Pi_1(y)$,
$y\in(0,x_*]$ enables us to apply
similar arguments to prove
that $\Delta(y)=0$,
$y\in(0,x_*)$.
\hfill$\square$

\subsection{Lower bounds}
This subsection provides a proof for the
existence of a strictly positive lower
bound for $\rho(x)$ in the $\ell$- and $s$-
cases. For $u$-demand, in  the particular
case $\sigma(s)= 1+\ln s$,
a weaker  estimate of
$\rho(x)=O((\ln x^{-1})^{-g})$, $g>0$ will be established
in the neighborhood of the origin (see
Lemma 4 in Appendix). This estimate
does not change the value of the Pareto
exponent, and because of Remark 1, the
necessary continuity condition $b_u=2$
also remains valid.

As
$\rho(x)$ is  continuous  and
$\rho(1)=1$, it follows that $\rho$ is
strictly positive  on a certain
segment $[x_*,1]$, i.e., there exists
$\rho_*$ such that
$$
\min_{x\in[x_*,1]}\rho(x)\ge\rho_*>0.
$$
For simplicity, suppose that   the point
$x_*$ is the same as that in equations
(\ref{estC}) and (\ref{MC}).
Let the $s$-demand be given by a smooth increasing
bounded function $\sigma(s)$ and that the tail of the
probability distribution $r(s)$ is not
too heavy, i.e., there exist $\mu$
such that
\begin{align}
\label{tail}
\frac{Z_\alpha(x)}{(1-x)\,c(x)}=
\frac{1}{(1-x)\,c(x)}\int_1^{1/x}s^\alpha r(s)\,ds\ge
\frac{1}{1+\mu\,x}\quad \forall
x\in(0,x_*].
\end{align}
For $\ell$- and $s$-demands, this estimate holds for
$r(s)$ with any tail dominated by $O(s^{-\alpha-2})$.
Under the same assumption on $r(s)$,
for $u$-demand we have
\begin{align}
\label{Utail}
\frac{Z_\alpha(x)}{(1-x)\,c(x)}\ge
\frac{1}{1+\mu\,x}\,\biggl(1+\frac{1}{R_0\sigma(x^{-1})}\biggr)^{-1}.
\end{align}
In the following theorem, the
lower bounds are established for $\rho$ under
assumptions (\ref{tail})--(\ref{Utail})
for any demand, but
in $u$-case, as  a typical demand
function
$\sigma(s)=1+\ln s$ will be
considered. The proof  remains valid for a
wider class of slowly varying functions that inherit
the property of the
logarithm $\sigma(\prod s_k)\le C\,\sum \sigma(s_k)$,
with a constant $C$ that is
independent of the number of factors
$\{s_k\}$.

\begin{theorem}
\label{LowerBound}
For $\ell$- and
$s$-demands,
$\rho(x)$ is strictly positive,  while for $u$-demand,
$\rho(x)$ is bounded below by $1/f(x^{-1})$, where
$f(s)$ a slowly varying function. More precisely,
$f(s)=(\ln s)^g$ and
$$
\rho(x)\ge \rho_*\,e^{-C}\,(\ln
x^{-1})^{-g}
$$
with positive constants $C$ and $g$
that are
specified in the proof.
\end{theorem}
\proof
The proof
does not differ for $\ell$- and $s$-cases.
Lemma 1 implies the
a.s. existence of a finite exit time
$n_*=n_*(x,x_*)$ from the segment
$(x,x_*]$ for the Markov chain (\ref{MC})
$\{x_n\}$ that starts at the point
$x_0=x$:  $x_{n_*-1}\le x_*<x_{n_*}$.
Suppose that $d\le0$.
Hence, we have
\begin{align*}
\rho(x)=&\e_x \prod_n F(x_{n-1},x_n)=
\e_x \e_{x_{n_*}}\prod_n
F(x_{n-1},x_n)
\\
&=\e_x \prod_{n=1}^{n_*-1}
F(x_{n-1},x_n)\;
\e_{x_{n_*-1}}F(x_{n_*-1},x_{n_*})\rho(x_{n_*})
\ge \rho_*\,
\e_x
\prod_{n=0}^{n_*-1}
\frac{1}{1+\mu
x_n}.
\end{align*}
Taking into account the inequalities
$\e\, e^\xi\ge e^{\e\,\xi}$ and $\ln(1+ x)\le
x$,
we obtain the estimate
\begin{align}
\label{EstExp} \rho(x)\ge
\rho_*\exp\biggl\{-\mu\;\e_x
\,\sum_{n=0}^{n_{*}}\,x_n\biggr\},\quad
x\in(0,x_*],
\end{align}
where  the points $x_n$ can be
expressed as the product of factors
$\xi_k=x_{k+1}/x_k\in (1,\infty)$, whose distribution
is given by the first equation
(\ref{MC}):
\begin{align*}
x_n=x\prod_{k=1}^{n}\frac{x_k}{x_{k-1}}=&
x_{n_*-1}\prod_{k=n-1}^{n_*-1}\biggl(\frac{x_k}{x_{k-1}}
\biggr)^{-1}\le
x_*\prod_{k=n-1}^{n_*-1}\xi_k^{-1},
\\
\p_x\{\xi>1/x\}=&\frac{1}{Z_\alpha(x)}
\int_1^{1/x}\xi^\alpha\,r(\xi)\,d\xi.
\end{align*}
Note that the function $\e_{x}\xi^{-1}$
increases in $x\in(0,x_*]$, i.e., $\e_{x}\xi^{-1}\le
\e_{x_*}\xi^{-1}\defn e_*<1$,  and
the Markov property implies the
estimate
$$
\e_x x_n=\e_x
x_{n_*-1}\prod_{k=n-1}^{n_*-1}\e_{x_k}\xi^{-1}\le
x_*\,\sum_n e_*^n\le (1-e_*)^{-1},
\quad x_*< 1.
$$
Hence, the series $\e_x\sum_n\,x_n$
converges uniformly in $x$; this
fact implies  the  strict positivity of
$\rho$:
$\rho(x)\ge \rho_*
\,\exp\{-\mu/(1-e_*)\}>0$, where $\mu $
is defined by (\ref{tail}).

Consider the case $d>0$. By the Markov
property, we have
\begin{align}
\label{d>0}
\e_x
\left(\frac{1-x_{n_*}}{1-x}\right)^d\;
\prod_{n=0}^{n_*-1}
\frac{1}{1+\mu
x_n}=\e_x
\prod_{n=0}^{n_*-1}
\frac{1}{1+\mu
x_n}\;\e_{x_{n_*-1}}\left(\frac{1-x_{n_*}}{1-x}\right)^d
\nonumber
\\
\ge
\inf_{x_m\le x_*}
\e_{x_m}(1-x_{m+1})^d\;
\e_x
\prod_{n=0}^{n_*-1}
\frac{1}{1+\mu
x_n},
\end{align}
where the expectation of the last
product is finite (see the proof for
$d\le0$) and the first infimum is
strictly positive:
\begin{align*}
\inf_{x_m\le x_*}
\e_{x_m}(1-x_{m+1})^d=\inf_{x\le x_*}
\frac{1}{Z_\alpha(x)}\int_1^{1/x}s^\alpha
r(s)\,(1-sx)^d\,ds
\\
\ge
\inf_{x\le x_*}
\frac{1}{Z_\alpha(x)}\frac{1}{2^d}\int_1^{1/2x}s^\alpha
r(s)\,ds=
\frac{Z_\alpha(2x_*)}{2^d Z_\alpha(x_*)}>0
\end{align*}
because the integral  normalized by $Z_\alpha$
in the last line of the equation is a
decreasing function of $x$.

Consider the $u$-case with
$\sigma(s)=1+\ln s$. For simplicity,
suppose that $d\le 0$; in the case
$d>0$, the proof is similar to that in
the previous paragraph. Under
the assumption $r(s)\le O(s^{-2})$, there
exists $\mu>0$ such that
$$
\frac{Z_\alpha(x)}{(1-x)c(x)}\ge\frac{1}{1+\mu x}
\,\frac{1}{1+1/R_0\ln x^{-1}}
$$
(see (\ref{Utail})).
By Lemma \ref{exitL}, the exit time can
be restricted by $N(x_*,x)$  with
a probability greater than or equal to
$(2Y^*(x)/Y_*-1)^{-1}$,
where $Y^*(x)=\ln(x_*/x)$ and
$Y_*=\e_{x_*}\ln\xi$.
Then,
similarly to (\ref{EstExp}), we get
\begin{align*}
\rho(x)\ge& \rho_*\e_x \exp\biggl\{-
 \sum_{n=0}^{n_*}\biggl(\mu x_n+
\frac{1}{R_0\ln x_n^{-1}}\biggr)\biggr\}
\\
&\ge
\frac{\rho_*}{2Y^*(x)/Y_*-1}\,
\e_x\exp\biggl\{-
  \sum_{n=0}^{N(x,x_*)}\biggl(\mu\, x_n+\frac{1}{R_0\,\ln x_n^{-1}}
  \biggr)\biggr\}
\\
&\ge
\frac{\rho_*}{2Y^*(x)/Y_*-1}\,
\exp\biggl\{-
  \sum_{n=0}^{N(x,x_*)}\biggl(\mu\,\e_x x_n+
  \e_x\frac{1}{R_0\,\ln x_n^{-1}}
  \biggr)\biggr\},
\end{align*}
where the first sum is uniformly
bounded in $x$ (see the above proof for the $\ell$- and $s$-cases).
As $\ln (1/x)$ is a concave function
on $x\in(0,e^{-2}]$, it is assumed that $x_*\le
e^{-2}$ and apply the Jensen
inequality to estimate the
last sum:
$$
\e_x\frac{1}{\ln x_n^{-1}}\le
\biggl(\ln \frac{1}{\e_x x_n}\biggr)^{-1}=
\biggl(\ln \frac{1}{x_{n_*}\e_x
\prod_{k=n}^N\xi_k^{-1}}\biggr)^{-1},\quad \xi_k=x_k/x_{k-1},
\;N=N(x,x_*).
$$
Note that the function $\e_x \xi^{-1}=\frac{1}{Z_
\alpha(x)}\int_{1}^{1/x}\xi^{-1+\alpha}\,r(s)\,ds$
increases in $x$. Hence, $\e_x
\xi^{-1}\le \e_{x_*} \xi^{-1}$, and
therefore
$$
\ln\frac{1}{x_{n_*}\e_x
\prod_{k=n}^N\xi_k^{-1}}\le\ln \frac{1}{x_{n_*}}+
\ln \frac{1}{\prod_{k=n}^{N(x,x_*)}\e_{x_*}
\xi^{-1}}\le
(N(x,x_*)-n+1)\,\lambda_*.
$$
where
$\lambda_*=\ln(\e_{x_*}\xi^{-1})^{-1}>0$.
This estimate enables completion of
the proof:
$$
 \sum_{n=0}^{N(x,x_*)}
  \e_x\frac{1}{R_0\,\ln x_n^{-1}}\le
\frac{1}{\lambda_*\,R_0}\sum_{n=1}^{N(x,x_*)}\frac{1}{n}\le
\frac{1+\ln
N(x,x_*)}{\lambda_*\,R_0},
$$
where $N(x,x_*)\le\frac{2}{Y_*}\ln
\frac{1}{x}$, i.e., the assertion of the
theorem
$\rho(x)\ge \rho_*\,e^{-C}\,(\ln
x^{-1})^{-g}$
holds with
$$
g=\frac{1}{\lambda_*\,R_0},\quad C=\frac{\mu}{1-e_*}+
\frac{1+\ln \frac{2}{Y_*}}{\lambda_*\,R_0}.
$$
It is clear that the factor $(\ln
x^{-1})^{-g}$, which
is a slowly varying function, does
not disturb the Pareto exponent.
 \hfill $\square$

\section{Conclusions}
The aim of this paper is to develop a
qualitative  model that implies stable
income distributions with credible
Pareto tails. To this end, the extreme
case of economical cooperation has been
considered and the basic  features of
chain-trading structures with rigid
corporative hierarchy have been
axiomatized. An analysis of
Section 2 shows that  the  qualitative
definition of  demand functions
determines the values of  Pareto
exponents for gross and net income. In
general, the choice of the demand
function is predetermined by the
economic environment and the scale of
the business. Hence, the qualitative behavior of
the demand function  can be motivated by
estimates of market size and other
macroeconomic constraints. The
observations that are valid in terms of
our model may provide some hints for
the interpretation of facts established
in empirical  studies. For example,
according to a number of reports (see
\cite{ASN00}--\cite{FuGuA05},
\cite{MKTT}--\cite{KII}), the graph of
the net income distribution drops down
from the Pareto exponent $a=1$ to a
greater value $1+\alpha$. This
phenomena can be explained as a
transition from a long-range $u$-demand
to $s$-limitation of the market reached
by  the largest
corporations.

The existence of the peak
of the density distribution in the
neighborhood of the minimal income
($d<0$) depends on  the amplitude
$R(1)$  of the pair-potential
 and the intensity of
interactions $R_0=\int R(s)\,ds$;
the accumulation of the poverty cannot be
controlled by  the only amplitude $P(1)$  of the
social welfare.  To remove the
peak of poverty (to make $d\ge 0$) or
at least reduce it, according to
equation (\ref{d}) it is better to
decrease the gradient of the social aid
($P'(1)\to 0)$, increase the
amplitude of pair-interactions $R(1)$
(i.e., to increase the economic
activity of agents), and to decrease
$R_0/R(1)=\int R(x)/R(1)\,dx$ (to
decrease the range of spending). These
expectable issues mean that the
mathematical caricature
(\ref{stat})--(\ref{minimal}) is
potentially consistent. The quotation ``{\it If I
give out my money to the poor in the
morning, all money will be repaid by
tonight\/}''  (apparently made by George
Soros) aphoristically expresses the
concept of the hierarchical model and motivates the basic
equations of this paper.

\section{Mathematical Appendix: Four lemmas
on multiplicative Markov chains}

\begin{lemma}
\label{exit} {\rm(On the a.s. finiteness of the exit time)}
The Markov chain $\{x_n\}$ with transition
probabilities {\rm(\ref{MC})} starting at
the point $x_0=x$ exits  any interval
$(x,y)\subset(0,1)$ in a finite number of
steps a.s.
\end{lemma}

\begin{proof}
(a) If $x\in(0,x_*]$, then by definition
(\ref{MC}),
$$
\p\{x_{n+1}\ge(1+\varepsilon)x_n|x_n\}=\frac{1}{x_nZ_\alpha(x_n)}
\int_{(1+\varepsilon)x_n}^{1}s^\alpha r(y/x)\,dy=
\frac{1}{Z_\alpha(x_n)}
\int_{1+\varepsilon}^{1/x_n}s^\alpha r(s)\,ds,
$$
where the right-hand side is a decreasing function of $x_n$,
as the derivative in $x_n$  of the
right-hand side of the equation is nonpositive.
If $\int_{1}^{1/x_*}r(s)\,ds=0$, the
Markov chain exits $(x,x_*]$  in the first step, with
a probability of one.
Otherwise, set
$\xi_{n+1}= x_{n+1}/x_n\ge 1$. As $\int_{1}^{1/x_*}
s^\alpha r(s)\,ds>0$
and $r\in L_1$,
there exists $\varepsilon>0$ such that
$$
\p\{\xi_{n+1}\ge 1+\varepsilon\}\ge
\frac{1}{Z_\alpha(x_*)}
\int_{1+\varepsilon}^{1/x_*}s^\alpha r(s)\,ds\ge
\delta_\varepsilon>0, \quad
\varepsilon< x_*^{-1}-1
$$
independently of the previous events.
The probability $P_\varepsilon(n,N)$ of the event
``{$n$ random variables $\xi_k$ among $N$
take values larger  than $1+\varepsilon$}''
can be estimated from the binomial
distribution:
\begin{align*}
P_\varepsilon(n,N)=&\p\bigl\{\prod_1^N\xi_k\ge(1+\varepsilon)^n\bigr\}\ge
\sum_{k=n+1}^N
C_N^k\,\delta_\varepsilon^k(1-\delta_\varepsilon)^{N-k}
\\
&=
1-\sum_{k=0}^n
C_N^k\,\delta_\varepsilon^k(1-\delta_\varepsilon)^{N-k}=
1-C_N^n\,\delta_\varepsilon^n(1-\delta_\varepsilon)^{N-n}
\frac{\delta_\varepsilon\,(N+1-n)}{(N+1)\delta_\varepsilon-n}
\\
&=
1-O\bigr(\exp\{n\ln N-N\ln
(1-\delta_\varepsilon)^{-1}\}\bigl)\to 1
\end{align*}
as $N\to\infty$ for any finite $n$
and $x<x_*$(\cite{Fel}, Ch.6, \S3).
Hence, the Markov chain exits the
interval
$(x,x_*]$ with probability one.
\medskip

(b) If $x\in(x_*,y]$, then by
(\ref{MC}), for  $\varepsilon<2/(1+y)$ we have
\begin{align*}
\p\{x_{n+1}\ge(1+\varepsilon)x_n|x_n\}=&
-\int_{(1+\varepsilon)x_n}^{1}\,d\left(\frac{1-y}{1-x}\right)^\gamma
=\left(1+\varepsilon-\frac{\varepsilon}{1-x_n}\right)^\gamma
\end{align*}
which is clearly a decreasing function of
$x_n$. Hence,
\begin{align*}
\p\{x_{n+1}\ge(1+\varepsilon)x_n|x_n\}\ge
\p\{x\ge(1+\varepsilon)y|y\}=
\left(1+\varepsilon-\frac{\varepsilon}{1-y}
\right)^\gamma=\delta_\varepsilon>0
\end{align*}
independently of the previous events.
By the same arguments as in case
(a),  for any finite $n$ and $x\ge x_*$, we have
$\p\bigl\{\prod_1^N\xi_k\ge(1+\varepsilon)^n|x\bigr\}\to 1$
as $N\to\infty$.
\hfill$\square$
\end{proof}

\begin{lemma} {\rm(On the finiteness of exponential
moments)}
\label{main}
Let $\{\xi_n\}$ be  i.i.r.v. with
moments $\mu_n=\e\,\xi^n<1$ and let $G$
be  the  random geometric progression
$$
G=\sum_{n=1}^\infty\prod_{k=1}^n
\,\xi_k=\xi_1\,(1+\xi_2\,(1+\xi_3\, (1+\,\dots))).
$$
Then
$
\e\,e^{\lambda G}<\infty$ for all
$\lambda\in \c$ such that
\begin{align}
\label{suff1}
\lambda: \quad\sigma(\lambda)\defn \sup_{k\ge1}\sum_{n=1}^\infty
\frac{|\lambda|^n}{n!}\,\frac{1}{1-\mu_{n+k}}<1.
\end{align}
\end{lemma}

\begin{proof} It suffices to prove the
statement in the case $\lambda\ge 0$.
Let us regard sum (\ref{geom}) as a
product of two  random
variables, $\xi$ and $1+\widetilde G$,
$$
G=
\xi_1\,(1+\xi_2\,(1+\dots))=
\xi_1\,(1+\widetilde G),
$$
where $\xi_1$  and $\widetilde G$ are
independent, $G$ and $\widetilde G$
are identically distributed. Denote by
$M_n$ the moments of the random
variable $G$: $M_n\defn \int
G^n\,P(dG)$. Therefore, the expansion
$$
\e\, e^{\lambda\,G}=\sum_{n=0}^{\infty}\frac{\lambda^n
\,M_n}{n!}=\int_\r \,\e\,
e^{\lambda\,\xi\,(1+\widetilde G)}\, d\xi=\sum_{n=0}^{\infty}
\frac{\lambda^n\,\mu_n
\,\e\,(1+G)^n}{n!}
$$
implies  a solvable  recurrent system of equations for the
moments $M_n$:
\begin{align}
\label{recurr}
M_n=\frac{1}{1-\mu_n}\,\sum_{k=0}^{n-1}C_n^k\,M_k,
\quad M_0=1,
\end{align}
 where $C_n^k=\frac{n!}{k!(n-k)!}$. Using
the recurrent definition (\ref{recurr})
of $M_n$, we obtain
\begin{align*}
S\defn\sum_{n=1}^{\infty}\frac{M_n\,\lambda^n}{n!}=&
\sum_{n=1}^{\infty}
\frac{\lambda^n}{n!}
\frac{1}{1-\mu_n}\biggl(1+\sum_{k=1}^{n-1}
M_k C_n^k\biggr)
\\
&=
\sum_{n=1}^{\infty}\frac{\lambda^n}{n!}\frac{1}{1-\mu_n}+
\sum_{k=1}^\infty\frac{\lambda^k\,M_k}{k!}\sum_{n=k+1}^\infty
\frac{\lambda^{n-k}}{(n-k)!}
\frac{1}{1-\mu_n}
\\
&=
\sum_{n=1}^{\infty}\frac{\lambda^n}{n!}\frac{1}{1-\mu_n}+
\sum_{k=1}^\infty\frac{\lambda^k\,M_k}{k!}\sum_{m=1}^\infty
\frac{\lambda^{m}}{m!}
\frac{1}{1-\mu_{k+m}}
\\
&\le L + S\, \sup_{k\ge 1}\sum_{m=1}^\infty
\frac{\lambda^{m}}{m!}
\frac{1}{1-\mu_{k+m}}= L+\sigma(\lambda) S,\quad
L\defn\sum_{n=1}^\infty\frac{\lambda^n}{n!\,(1-\mu_n)}.
\end{align*}
Thus, $S\le
L\,(1-\sigma(\lambda))^{-1}<\infty$,
provided $\sigma(\lambda)<1$.
\hfill$\square$
\end{proof}

\begin{lemma} {\rm(Anti-Chebyshev inequality)}
\label{Cheby}
Let $\{Y_n\}$ be a sequence of positive
finite random variables such that $Y_n\le Y^*<\infty$ and
$\inf_{n,Y_1,\dots,Y_m}\e(Y_{n+1}|Y_1,\dots,Y_m)=Y_*>0$.
Then,
$$
\p\biggl\{\sum_{k=1}^n Y_k>\frac{1}{2}n\,Y_*\biggr\}\ge
\frac{1}{2Y^*/Y_*-1}.
$$
\end{lemma}
\proof For $X=\sum_{k=1}^n Y_k$, we
have $n\,Y^*\ge X^*\ge \e\, X\ge n\,Y_*$.
Hence, for any $\varepsilon\in(0,1)$, we
have  $\p\{X>\varepsilon\, n\,Y_*\}\ge
\p\{X>\varepsilon\,\e X\}$. In addition,
$$
\e X\le X^*\,\p\{X>\varepsilon\,\e X\}+\varepsilon\,\e
X\,(1-\p\{X>\varepsilon\,\e X\}).
$$
Therefore, by setting $\varepsilon=1/2$,
we obtain
$$
\p\{X>\varepsilon\, n\,Y_*\}\ge
\p\{X>\varepsilon\,\e X\}\ge \frac{1-\varepsilon}{X^*/\e X-\varepsilon}
\ge
\frac{1-\varepsilon}{Y^*/Y_*-\varepsilon}\biggl|_{\varepsilon=1/2}
=\frac{1}{2Y^*/Y_*-1}.\qquad\square
$$

\begin{lemma}{\rm (On the distribution of the exit time)}
\label{exitL}
The Markov chain $\{x_n\}$ with transition probabilities
{\rm(\ref{MC})} exits the interval
$(x,x_*)$ in $N(x,x_*)$ steps with the probability
$$
\p\{n_*(x)\ge
N(x,x_*)\}\ge \frac{1}{2Y^*(x)/Y_*-1}=O(\ln
x^{-1})^{-1},
$$
where $Y^*(x)=\ln (x_*/x)\le\ln
x^{-1}$, and
$$
N(x,x_*)=\biggl[\frac{2}{Y_*}\ln
\frac{x_*}{x}\biggr]+1=O(\ln
x^{-1}),\quad
Y_*=\frac{1}{Z_\alpha(x_*)}\int_1^{1/x^*}\ln
s\,s^\alpha\,r(s)\,ds>0.
$$
\end{lemma}
\proof
Note that the function
$$
\e_{x_n}\ln
\frac{x_{n+1}}{x_n}=\frac{1}{Z_\alpha(x_n)}\int_1^{1/x_n}\ln
\xi\;\xi^\alpha\,r(\xi)\,d\xi\defn\e_{x_n}\ln
\xi
$$
is a decreasing function of the initial
point $x_n\in(0,x_*]$; hence, $\e_{x_n}\ln
\frac{x_{n+1}}{x_n}\ge \e_{x_*}\ln
\xi$, and for the random variables $Y_n=
\ln
\frac{x_{n+1}}{x_n}$, we have
$$
Y_*=\min\e_x\{Y_{n+1}|Y_1,\dots,Y_n\}=\e_{x_*}\ln
\xi.
$$
Hence, for  $x_n=x\,\prod_{k=1}^n\,x_k/x_{k-1}=x\,\exp\{
\sum_{k=0}^{n-1}\ln Y_k\}$,  Lemma
\ref{Cheby} implies
\begin{align*}
\p\{x_{n_*}>x_*\}&=\p\biggl\{\sum_{k=0}^{n_*-1}\ln Y_k>\ln x_*/x
\biggr\}
\\
&\ge\p\biggl\{\sum_{k=0}^{n_*-1}\ln
Y_k>n_*\,Y_*/2\biggr\}\biggl|_{n_*Y_*>2\ln (x_*/x)}
=\frac{1}{2Y^*(x)/Y_*-1},
\end{align*}
i.e., exit time $N(x,x_*)=\bigl[\frac{2}{Y_*}\ln
\frac{x_*}{x}\bigr]+1$ has a
probability of not less than
$(2Y^*(x)/Y_*-1)^{-1}$.

\hfill$\square$


\begin{thebibliography}{99}
\bibitem{Ch53}
D.~G.~Champernowne,
%A model of income distribution,
Econ. J., {\bf 23} (1953) 318--351.

\bibitem{MaZa99}
S.~C.~Manrubia, D.~H.~Zanette,
%Stochastic
%multiplicative processes with reset events.
Phys. Rew. E, {\bf 59} (1999)
4945--4948.

\bibitem{NS04}
M.~Nieri, W.~Suoma,
%Two-factors model of income
%distribution dynamics.
Working paper of
Dept. of Economics of Utah State
Univ., October 6 (2004).

\bibitem{Yul37}
G.~U.~Yule,
%A mathematical theory of evolution
%based on conclusions of Dr.
%J.~C.~Willis,
Philos. Trans. R. Soc. London B,
{\bf 213} (1925) 21--87.

\bibitem{New05} M.~E.~J. Newman,
%Power laws, Pareto distribution and Zipf's law
Contemporary Physics, {\bf 46}  (2005) 323--351.

\bibitem{AoSouF03}
H.~Aoyama, W.~Souma, Y.~Fujiwara,
%Growth and fluctuations of personal and company's income
%{\it Physica A} {\bf 324} (2003) 352--358.
 Physica A, {\bf 324} (2003)
352-358.
\bibitem{S39} M.~Y. Sweezy,
%Distribution of wealth and incomes
%under the Nazis, Rev.
Economic
Statistics, {\bf 21}, no. 4  (1939) 178--184.

\bibitem{ASN00}
H.~Aoyama, W.~Souma, Y.~Nagahara,
M.~P.~Okazaki, H.~Takayasu, M.~Takayasu,
%Pareto's law for income of individuals
%and debt of bankrupt companies,
Fractals,  {\bf 8},  no. 3 (2000) 293--300.

\bibitem{FuGuA05} Y.~Fujiwara, C.~Di Guilmi,
H.~Aoyama, M.~Galetti, W.~Souma,
%Do Pareto-Zipf and Gibrat laws hold true?
%An analysis with European Firms.
arXiv:cond-mat/0310061, 2003.

\bibitem{A04}
J.~H.~Ausubel,
%Will the rest of the world live like America?
Technology in Society, {\bf 26}
(2004) 343--360.

\bibitem{AtSal02}
A. B. Atkinson, % and Wiemer
W. Salverda,
%Top Incomes in the Netherlands and the United Kingdom over the Twentieth
%Century,
Journal of European Economic Assiciation,
{\bf 3}, no. 4  ( (2005)   883--913.

\bibitem{At03} A.~B.~Atkinson,
Income Tax and Top Incomes over the
Twentieth Century, Plenary Lecture
given at the XXVIII Meeting on Economic
Analysis, Seville, 11--12 December 2003.

\bibitem{R03}
W.~J.~Reed,
%The Pareto law of income---an
%explanation and an extension,
Physica A,
{\bf 319} (2003) 469--486.

\bibitem{Gab99}
X.~Gabaix,
%Zipf's law for cities: an explanation.
Quarterly Journal of Economics,
{\bf CXIV} (1999)
739--767.

%d=\infti, a>1
\bibitem{Bush01}
J.-Ph.~Bouchaud,
Power laws in economics and finance:
some ideas from physics,
Quantitative Finance, {\bf 1} (2001)
105--112.

\bibitem{KK03} %Christian
C.~Kleiber, %Samuel
S.~Kotz,
Statistical Size Distributions in Economics and Actuarial Sciences,
John Wiley \& Sons, New York, 2003.


\bibitem{MKTT}
T.~Mizuno, M.~Katori, H.~Takayasu,
M.~Takayasu,
%Statistical Laws in the Income of Japanese
%Companies,
In:
Empirical Science of Financial Fluctuations - The Advent of Econophysics,
Springer (2001) 321--330.

\bibitem{KII}
T.~Kaizoji,
Y.~Ikeda,
H.~Iyetomi,
%Re-examination of the size distribution of firms,
arXiv:physics/0512124, 2005.


\bibitem{Ch03}
A.~M.~Chebotarev,
%Pareto distribution as
%the result of computer reconstruction
%of Russian auto market statistics,
Ukrainian Journal ``Economist'',
no. 7  (2003) 6--11.

\bibitem{Ch05}
V.~N.~Baturin, S.~G.~Lebedev, V.~P.~Maslov,
B.~I.~Sadovnikov, A.~M.~Chebotarev,
%A conjecture on distribution of high
%incomes and its interpretation,
Economics of Contemporary Russia, no. 4,
(31) (2005) 57--62.

\bibitem{Ch04}
A.~M.~Chebotarev,
%Pareto distribution of incomes in
%hierarchical model of economy has
%exponent value of 2,
Ukrainian Journal ``Economist'',
no. 9 (2004) 54--57.

\bibitem{Bog99}
L.~Boggio,
%-Luciano, (I-SACM)
%On local relative stability of large systems with small parameters.
%The example of a classical model of competition.
Riv. Mat. Sci. Econom. Social.
%[Rivista-di-Matematica-per-le-Scienze-Economiche-e-Sociali]
22 (1999), no. 1-2, 13--30.


\bibitem{Gast72}
J.~L.~Gastwirth,
%The estimatiom of the
%Lorentz curve and Gini index,
Review of Economics and Statistics,
{\bf 54}, no. 3, (1972) 306--316.

\bibitem{Fel}
W.~Feller,
{\it An Introduction to Probability Theory and Its
Applications\/},
Vol.~1,
John Wiley \& Sons,
New York,
1966.

\end{thebibliography}
\end{document}